\documentclass{amsart}
\usepackage{amsopn}
\usepackage{amssymb, amscd}

\newcommand{\nc}{\newcommand}

\nc{\vg}{\mathfrak{v} } \nc{\wg}{\mathfrak{w} } \nc{\zg}{\mathfrak{z} } \nc{\ngo}{\mathfrak{n} }
\nc{\kg}{\mathfrak{k} } \nc{\mg}{\mathfrak{m} } \nc{\bg}{\mathfrak{b} } \nc{\ggo}{\mathfrak{g} }
\nc{\ggob}{\overline{\mathfrak{g}} } \nc{\sog}{\mathfrak{so} } \nc{\sug}{\mathfrak{su} } \nc{\spg}{\mathfrak{sp}
} \nc{\slg}{\mathfrak{sl} } \nc{\glg}{\mathfrak{gl} } \nc{\cg}{\mathfrak{c} } \nc{\rg}{\mathfrak{r} }
\nc{\hg}{\mathfrak{h} } \nc{\tg}{\mathfrak{t} } \nc{\ug}{\mathfrak{u} } \nc{\dg}{\mathfrak{d} }
\nc{\ag}{\mathfrak{a} } \nc{\pg}{\mathfrak{p} } \nc{\sg}{\mathfrak{s} } \nc{\pca}{\mathcal{P}}
\nc{\nca}{\mathcal{N}} \nc{\lca}{\mathcal{L}} \nc{\oca}{\mathcal{O}} \nc{\mca}{\mathcal{M}}
\nc{\tca}{\mathcal{T}} \nc{\aca}{\mathcal{A}} \nc{\cca}{\mathcal{C}} \nc{\sca}{\mathcal{S}}

\nc{\vp}{\varphi} \nc{\ddt}{{\small \frac{{\rm d}}{{\rm d}t}}} \nc{\im}{\mathtt{i}}

\nc{\SO}{{\mathrm SO}} \nc{\Spe}{{\mathrm Sp}} \nc{\Sl}{{\mathrm SL}} \nc{\SU}{{\mathrm SU}} \nc{\Or}{{\mathrm
O}} \nc{\U}{{\mathrm U}} \nc{\Gl}{{\mathrm GL}} \nc{\Se}{{\mathrm S}} \nc{\Cl}{{\mathrm Cl}}
\nc{\Spein}{{\mathrm Spin}} \nc{\Pin}{{\mathrm Pin}}

\nc{\RR}{{\Bbb R}} \nc{\HH}{{\Bbb H}} \nc{\CC}{{\Bbb C}} \nc{\ZZ}{{\Bbb Z}} \nc{\FF}{{\Bbb F}} \nc{\NN}{{\Bbb
N}} \nc{\QQ}{{\Bbb Q}} \nc{\PP}{{\Bbb P}}

\nc{\vs}{\vspace{.2cm}} \nc{\vsp}{\vspace{1cm}} \nc{\ip}{\langle\cdot,\cdot\rangle} \nc{\la}{\langle}
\nc{\ra}{\rangle} \nc{\unm}{\frac{1}{2}} \nc{\unc}{\frac{1}{4}} \nc{\und}{\frac{1}{16}} \nc{\no}{\vs\noindent}
\nc{\lam}{\Lambda^2\ngo^*\otimes\ngo} \nc{\tangz}{{\rm T}^{\rm Zar}} \nc{\nor}{{\sf n}}
\nc{\eigen}{(k_1<...<k_r;d_1,...,d_r)} \nc{\eigencero}{(0<k_2<...<k_r;d_1,...,d_r)} \nc{\mum}{/\!\!/}
\nc{\kir}{/\!\!/\!\!/}

\nc{\He}{\operatorname{Hess}} \nc{\ad}{\operatorname{ad}} \nc{\Ad}{\operatorname{Ad}}
\nc{\rank}{\operatorname{rank}} \nc{\Irr}{\operatorname{Irr}} \nc{\End}{\operatorname{End}}
\nc{\Aut}{\operatorname{Aut}} \nc{\Inn}{\operatorname{Inn}} \nc{\Der}{\operatorname{Der}}
\nc{\Ker}{\operatorname{Ker}} \nc{\Iso}{\operatorname{I}} \nc{\Diff}{\operatorname{D}}
\nc{\Lie}{\operatorname{L}} \nc{\tr}{\operatorname{tr}} \nc{\dif}{\operatorname{d}}
\nc{\sen}{\operatorname{sen}} \nc{\modu}{\operatorname{mod}} \nc{\Ric}{\operatorname{Ric}}
\nc{\Ricac}{\operatorname{Ric^{ac}}} \nc{\Ricg}{\operatorname{Ric^{\gamma}}} \nc{\Ricc}{\operatorname{Ric^{c}}}
\nc{\sym}{\operatorname{sym}} \nc{\symac}{\operatorname{sym^{ac}}} \nc{\symc}{\operatorname{sym^{c}}}
\nc{\scalar}{\operatorname{sc}} \nc{\grad}{\operatorname{grad}} \nc{\ricci}{\operatorname{ric}}
\nc{\ricciac}{\operatorname{ric^{ac}}} \nc{\riccic}{\operatorname{ric^{c}}}
\nc{\riccig}{\operatorname{ric^{\gamma}}} \nc{\Rin}{\operatorname{M}} \nc{\Le}{\operatorname{L}}
\nc{\tang}{\operatorname{T}} \nc{\level}{\operatorname{level}} \nc{\rad}{\operatorname{r}}
\nc{\abel}{\operatorname{ab}}

\newtheorem{theorem}{Theorem}[section]
\newtheorem{proposition}[theorem]{Proposition}
\newtheorem{corollary}[theorem]{Corollary}
\newtheorem{lemma}[theorem]{Lemma}
\newtheorem{definition}[theorem]{Definition}
\newtheorem{remark}[theorem]{Remark}

\newtheorem{example}[theorem]{Example}

\title[Geometric structures on nilpotent Lie groups]{Geometric structures on nilpotent Lie groups: on
their classification and a distinguished compatible metric}

\author{Jorge Lauret}

\address{Department of Mathematics, Yale University,
10 Hillhouse Box 208283 New Haven, CT 06520 USA (current affiliation)} \email{jorge.lauret@yale.edu}

\address{FaMAF and CIEM, Universidad Nacional de C\'ordoba, 5000 C\'ordoba, Argentina}
\email{lauret@mate.uncor.edu}

\thanks{2000 {\it Mathematics Subject Classification.} Primary: 53D05, 53D55;
Secondary: 22E25, 53D20, 14L24, 53C30. \\
{\it Key words and phrases.}  symplectic, complex, hypercomplex, nilpotent Lie groups,
moment map, variety of Lie algebras. \\
Supported by CONICET and Guggenheim Foundation fellowships, and a grant from FONCyT (Argentina).}

\begin{document}

\maketitle

\begin{abstract}
Let $(N,\gamma)$ be a nilpotent Lie group endowed with an invariant geometric structure (cf. symplectic,
complex, hypercomplex or any of their `almost' versions).  We define a left invariant Riemannian metric on $N$
compatible with $\gamma$ to be {\it minimal}, if it minimizes the norm of the invariant part of the Ricci tensor
among all compatible metrics with the same scalar curvature.  We prove that minimal metrics (if any) are unique
up to isometry and scaling, they develop soliton solutions for the `invariant Ricci' flow and are characterized
as the critical points of a natural variational problem. The uniqueness allows us to distinguish two geometric
structures with Riemannian data, giving rise to a great deal of invariants.

Our approach proposes to vary Lie brackets rather than inner products; our tool is the moment map for the action
of a reductive Lie group on the algebraic variety of all Lie algebras, which we show to coincide in this setting
with the Ricci operator.  This gives us the possibility to use strong results from geometric invariant theory.
We describe the moduli space of all isomorphism classes of geometric structures on nilpotent Lie groups of a
given class and dimension admitting a minimal compatible metric, as the disjoint union of semi-algebraic
varieties which are homeomorphic to categorical quotients of suitable linear actions of reductive Lie groups.
Such special geometric structures can therefore be distinguished by using invariant polynomials.
\end{abstract}

\section{Introduction}\label{intro}

Invariant structures on nilpotent Lie groups, as well as on their compact versions, nilmanifolds, play an
important role in symplectic and complex geometry.  The aims of this paper are the search for the `best'
compatible metric and the classification of such structures up to isomorphism, which are, in the end, two
intimately related problems.  Symplectic, complex and hypercomplex cases, and their respective `almost'
versions, will be treated in some detail, and although we will always have these particular cases in mind, the
main results will be proved for a geometric structure in general.  Contact and complex symplectic structures
will be studied in a forthcoming paper.

\subsection{Geometric structures and compatible metrics}
Let $N$ be a real $n$-di\-men\-sio\-nal nilpotent Lie group with Lie algebra $\ngo$, whose Lie bracket will be
denoted by $\mu :\ngo\times\ngo\mapsto\ngo$.  An invariant geometric structure on $N$ is defined by left
translation of a tensor $\gamma$ on $\ngo$ (or a set of tensors), usually non-degenerate in some way, which
satisfies a suitable integrability condition
\begin{equation}\label{closedgint}
{\rm IC}(\gamma,\mu)=0,
\end{equation}
involving only $\mu$ and $\gamma$.  The pair $(N,\gamma)$ will be called a {\it class-$\gamma$ nilpotent
Lie group}, and $N$ will be assumed to be simply connected for simplicity.  A left invariant Riemannian metric
on $N$ is said to be {\it compatible} with $(N,\gamma)$ if the corresponding inner product $\ip$ on $\ngo$
satisfies an orthogonality condition
\begin{equation}\label{ortconint}
{\rm OC}(\gamma,\ip)=0,
\end{equation}
in which only $\ip$ and $\gamma$ are involved.  We denote by $\cca=\cca(N,\gamma)$ the set of all left invariant
metrics on $N$ which are compatible with $(N,\gamma)$.  The pair $(\gamma,\ip)$ with $\ip\in\cca$ will often be
referred to as a {\it class-$\gamma$ metric structure}.

A natural question arises:
\begin{itemize}
\item[] Given a class-$\gamma$ nilpotent Lie group $(N,\gamma)$, are there distinguished left invariant
Riemannian metrics on $N$ compatible with $\gamma$?,
\end{itemize}
where the meaning of `distinguished' is of course part of the problem.  The Ricci tensor has always been a very
useful tool to deal with this kind of questions, and since the answer should depend on the metric and on the
structure under consideration, we consider the {\it invariant Ricci operator} $\Ricg_{\ip}$ (and the invariant
Ricci tensor $\riccig_{\ip}=\la\Ricg_{\ip}\cdot,\cdot\ra$), that is, the orthogonal projection of the Ricci
operator $\Ric_{\ip}$ onto the subspace of those symmetric maps of $\ngo$ leaving $\gamma$ invariant. D. Blair,
S. Ianus and A. Ledger \cite{BlrIns, BlrLdg} have proved in the compact case that metrics satisfying
\begin{equation}\label{unconditionint}
\riccig_{\ip}=0
\end{equation}
are very special in symplectic (so called metrics with hermitian Ricci tensor) and contact geometry, as they are
precisely the critical points of two very natural curvature functionals on $\cca$: the total scalar curvature
functional $S$ and a functional $K$ measuring how far are the metrics of being K$\ddot{{\rm a}}$hler or
Sasakian, respectively (see also \cite{Blr} and Section \ref{hermricci}).

We will show that for a non-abelian nilpotent Lie group, condition (\ref{unconditionint}) cannot hold for the
classes of structures we have in mind, and hence it is natural to try to get as close as possible to this
unnatainable goal.  In this light, a metric $\ip\in\cca(N,\gamma)$ is called {\it minimal} if it minimizes the
functional $||\riccig_{\ip}||^2=\tr(\Ricg_{\ip})^2$ on the set of all compatible metrics with the same scalar
curvature.  It turns out that minimal metrics are the elements in $\cca$ closest to satisfy the `Einstein-like'
condition $\riccig_{\ip}=c\ip$, $c\in\RR$.  We may also try to improve the metric via the evolution flow
$$
\ddt \ip_t=\pm\riccig_{\ip_t},
$$
whose fixed points are precisely metrics satisfying (\ref{unconditionint}).  In the symplectic case, this flow
is called the anticomplexified Ricci flow and has been recently studied by H-V Le and G. Wang \cite{LeWng}.  Of
particular significance are then those metrics for which the solution to the normalized flow (under which the
scalar curvature is constant in time) remains isometric to the initial metric.  Such special metrics will be
called {\it invariant Ricci solitons}.  The main result in this paper can be now stated.

\begin{theorem}\label{equiv2gint}
Let $(N,\gamma)$ be a nilpotent Lie group endowed with an invariant geometric structure $\gamma$
(non-necessarily integrable). Then the following conditions on a left invariant Riemannian metric $\ip$ which is
compatible with $(N,\gamma)$ are equivalent:
\begin{itemize}
\item[(i)] $\ip$ is minimal.

\item[(ii)] $\ip$ is an invariant Ricci soliton.

\item[(iii)] $\Ricg_{\ip}=cI+D$ for some $c\in\RR$, $D\in\Der(\ngo)$.
\end{itemize}
Moreover, there is at most one compatible left invariant metric on $(N,\gamma)$ up to isometry  (and scaling)
satisfying any of the above conditions.
\end{theorem}

A major obstacle to classify geometric structures is the terrible lack of invariants.  The uniqueness result in
the above theorem gives rise to a very useful tool to distinguish two geometric structures; indeed, if they are
isomorphic then their respective minimal compatible metrics (if any) have to be isometric.  One therefore can
eventually distinguish geometric structures with Riemannian data, which suddenly provides us with a great deal
of invariants.  This will be used to find explicit continuous families depending on $1$, $2$ and $3$ parameters
of pairwise non-isomorphic geometric structures in low dimensions, mainly by using only one Riemannian
invariant: the eigenvalues of the Ricci operator (see Sections \ref{sint}-\ref{hint} for a summary in each
case).

A remarkable weakness of this approach is however the existence problem; the theorem does not even suggest when
such a distinguished metric does exist.  How special are the symplectic or (almost-) complex structures
admitting a minimal metric?.  So far, we know how to deal with this `existence question' only by giving several
explicit examples, for which the neat `algebraic' characterization (iii) will be very useful. It turns out that
in low dimensions the structures in general tend to admit a minimal compatible metric, and the only obstruction
we know at this moment is when $G_{\gamma}\subset\Sl(n)$ and the nilpotent Lie algebra is characteristically
nilpotent, that is, $\Der(\ngo)$ is nilpotent.  Anyway, we could not find a non-existence example yet.

\subsection{Variety of compatible metrics and the moment map}\label{varint}
A class-$\gamma$ metric structure on a nilpotent Lie group is determined by a triple $(\mu,\gamma,\ip)$ of
tensors on $\ngo$.  The proof of Theorem \ref{equiv2gint} is based on an approach which proposes to vary the Lie
bracket $\mu$ rather than the inner product $\ip$.

Let us consider as a parameter space for the set of all real nilpotent Lie algebras of a given dimension $n$,
the set
$$
\nca=\{\mu\in V:\mu\;\mbox{satisfies Jacobi and is nilpotent}\},
$$
where $\ngo$ is a fixed $n$-dimensional real vector space and $V=\lam$ is the vector space of all skew-symmetric
bilinear maps from $\ngo\times\ngo$ to $\ngo$. Since the Jacobi identity and the nilpotency condition are both
determined by zeroes of polynomials, $\nca$ is an algebraic subset of $V$.  We fix a tensor $\gamma$ on $\ngo$
(or a set of tensors), and denote by $G_{\gamma}$ the subgroup of $\Gl(n)$ preserving $\gamma$. The reductive
Lie group $G_{\gamma}$ acts naturally on $V$ leaving $\nca$ invariant and also the algebraic subset
$\nca_{\gamma}\subset\nca$ given by
$$
\nca_{\gamma}=\{\mu\in\nca:{\rm IC}(\gamma,\mu)=0\},
$$
that is, those nilpotent Lie brackets for which $\gamma$ is integrable (see (\ref{closedgint})).

For each $\mu\in\nca$, let $N_{\mu}$ denote the simply connected nilpotent Lie group with Lie algebra
$(\ngo,\mu)$.  Fix an inner product $\ip$ on $\ngo$ compatible with $\gamma$, that is, such that
(\ref{ortconint}) holds.  We identify each $\mu\in\nca_{\gamma}$ with a class-$\gamma$ metric structure on a
nilpotent Lie group
\begin{equation}\label{idegint}
\mu\longleftrightarrow   (N_{\mu},\gamma,\ip),
\end{equation}
where all the structures are defined by left invariant translation.  The orbit $G_{\gamma}.\mu$ parameterizes
then all the left invariant metrics which are compatible with $(N_{\mu},\gamma)$ and hence we may view
$\nca_{\gamma}$ as the space of all class-$\gamma$ metric structures on nilpotent Lie groups of dimension $n$.
Two metrics $\mu,\lambda\in\nca_{\gamma}$ are isometric if and only if they live in the same $K_{\gamma}$-orbit,
where $K_{\gamma}=G_{\gamma}\cap\Or(\ngo,\ip)$.

We now go back to our search for the best compatible metric.  It is natural to consider the functional
$F:\nca_{\gamma}\mapsto\RR$ given by $F(\mu)=\tr(\Ricg_{\mu})^2$, which in some sense measures how far the
metric $\mu$ is from satisfying (\ref{unconditionint}).  The critical points of $F/||\mu||^4$ on the projective
algebraic variety $\PP\nca_{\gamma}\subset\PP V$ (which is equivalent to normalize by the scalar curvature since
$\scalar(\mu)=-\unc ||\mu||^2$), may therefore be considered compatible metrics of particular significance.

A crucial fact of this approach is that the moment map $m_{\gamma}:V\mapsto\pg_{\gamma}$ for the action of
$G_{\gamma}$ on $V$, where $\pg_{\gamma}$ is the space of symmetric maps of $(\ngo,\ip)$ leaving $\gamma$
invariant (i.e. $\ggo_{\gamma}=\kg_{\gamma}\oplus\pg_{\gamma}$ is a Cartan decomposition), satisfies
$$
m_{\gamma}(\mu)=8\Ricg_{\mu}, \qquad \forall\;\mu\in\nca_{\gamma},
$$
where $\Ricg_{\mu}$ is the invariant Ricci operator of $\mu$.  This allows us to use strong and well-known
results on the moment map due to F. Kirwan \cite{Krw1} and L. Ness \cite{Nss}, and proved by A. Marian
\cite{Mrn} in the real case (see Section \ref{cpmp} for an overview).  Indeed, since $F$ becomes a scalar
multiple of the square norm of the moment map, we obtain the following

\begin{theorem}\label{equiv1gint}\cite{Mrn}
Let $F:\PP \nca_{\gamma}\mapsto\RR$ be defined by $F([\mu])=\tr(\Ricg_{\mu})^2/||\mu||^4$.  Then for $\mu\in
\nca_{\gamma}$ the following conditions are equivalent:
\begin{itemize}
\item[(i)] $[\mu]$ is a critical point of $F$.

\item[(ii)] $F|_{G_{\gamma}.[\mu]}$ attains its minimum value at $[\mu]$.

\item[(iii)] $\Ricg_{\mu}=cI+D$ for some $c\in\RR$, $D\in\Der(\mu)$.
\end{itemize}
Moreover, all the other critical points of $F$ in the orbit $G_{\gamma}.[\mu]$ lie in $K_{\gamma}.[\mu]$.
\end{theorem}

Theorem \ref{equiv2gint} follows then almost directly from this result, except for the equivalence between (ii)
and (iii), which will be proved separately.  We note that Theorem \ref{equiv1gint} also gives a variational method to
find minimal compatible metrics, by characterizing them as the critical points of a natural curvature functional
(see Example \ref{m26} for an explicit application).

Most of the results obtained in this paper are still valid for general Lie groups, although some considerations
have to be carefully taken into account (see Remark \ref{liegroups}).

\subsection{Symplectic structures}\label{sint}
We first prove that a symplectic non-abelian nilpotent Lie group $(N,\omega)$ can never admit a compatible left
invariant metric with hermitian Ricci tensor.  We also find a minimal compatible metric for the two
$4$-dimensional symplectic nilpotent Lie groups and exhibit curves of pairwise non-isomorphic symplectic
structures on the $6$-dimensional nilpotent Lie groups denoted by $(0,0,0,12,13,23)$ and
$(0,0,12,13,14+23,24+15)$ in \cite{Slm}.

\subsection{Complex structures}\label{cint}
We find two different curves of pairwise non-isomorphic abelian complex structures on the Iwasawa manifold; for
only one of them their minimal compatible metrics are modified H-type. A third curve of pairwise non-isomorphic
non-abelian complex structures on the Iwasawa manifold is also given.  The initial point is the bi-invariant
complex structure, which is the only point for which the minimal metric is modified H-type. We also exhibit a
curve of pairwise non-isomorphic abelian complex structures on $\hg_3\oplus\hg_3$, where $\hg_3$ is the
$3$-dimensional Heisenberg Lie algebra, and a curve of non-abelian ones on the group denoted by
$(0,0,0,0,12,14+23)$ in \cite{Slm}.

\subsection{Hypercomplex structures}\label{hint}
We prove that any hypercomplex $8$-dimensional nilpotent Lie group admits a minimal compatible metric (being
actually the only compatible metric up to isometry and scaling), which is modified H-type if and only if the
structure is in addition abelian.  Let $\ggo_1$, $\ggo_2$ and $\ggo_3$ denote the $8$-dimensional Lie algebras
obtained as the direct sum of an abelian factor and the $5$-dimensional Heisenberg Lie algebra, the
$6$-dimensional complex Heisenberg Lie algebra and the $7$-dimensional quaternionic Heisenberg Lie algebra,
respectively.  A curve and a surface of pairwise non-isomorphic abelian hypercomplex structures on $\ggo_2$ and
$\ggo_3$ respectively, are given.  We also find curves of pairwise non-isomorphic non-abelian hypercomplex
structures on $\ggo_3$ and $\ngo=\ug(2)\oplus\CC^2$.  By using results due to I. Dotti and A. Fino
\cite{DttFin0, DttFin2}, we actually prove that the moduli space of all hypercomplex $8$-dimensional nilpotent
Lie groups, up to isomorphism, is parameterized by the $9$-dimensional quotient
$$
\PP\left( (\sug(2)\otimes\RR^4)\oplus\RR^4\right)/(\SU(2)\times\SU(2)),
$$
where $\sug(2)$ is the adjoint representation and $\RR^4$ is the standard representation of $\SU(2)$ on $\CC^2$
viewed as real.  The abelian ones are parameterized by the quotient
$$
\PP\left( \sug(2)\otimes\RR^4\right)/(\SU(2)\times\SU(2)),
$$
which has dimension $5$.  Explicit continuous families depending on $5$ parameters on $\ggo_2$ and $\ggo_3$ are
given.

\subsection{Einstein solvmanifolds}\label{eint}
If one considers no structure (i.e. $\gamma=0$), then we show that the `moment map' approach proposed in this
paper can be also applied to the study of Einstein solvmanifolds.  Each $\mu\in\nca$ is identified via
(\ref{idegint}) with the Riemannian manifold $(N_{\mu},\ip)$, but we also have in this case another
identification with a solvmanifold: for each $\mu\in\nca$, there exists a unique rank-one metric solvable
extension $S_{\mu}=(S_{\mu},\ip)$ of $(N_{\mu},\ip)$ standing a chance of being Einstein, and
every $(n+1)$-dimensional rank-one Einstein solvmanifold can be modelled as $S_{\mu}$ for a suitable
$\mu\in\nca$.  The functional $F$ measures how far is the metric $\mu$ from being Einstein.  We obtain, as a
consequence of the above theorems, many of the uniqueness and structure results proved by J. Heber in
\cite{Hbr}, as well as the variational result in \cite{critical} and the relationship between Ricci soliton
metrics on nilpotent Lie groups and Einstein solvmanifolds proved in \cite{soliton}.

\subsection{On the classification of geometric structures}\label{clint}
Everything seems to indicate that the moduli space of isomorphism classes of $n$-dimensional, class-$\gamma$,
nilpotent Lie groups is a very complicated space for most of the classes of geometric structures, even in low
dimensions.  Anyway, what can be said about such a moduli space?. Can we show that it is really unmanageable?.
Can we at least find subspaces which are manifolds or algebraic varieties and obtain lower bounds for its
`dimension'?. This kind of questions belong to invariant theory. For a fixed nilpotent Lie group $N$, the
isomorphism between geometric structures is determined by the action of $\Aut(\ngo)$, which is a group in
general unknown and `very ugly' from an invariant-theoretic point of view since it is far from being semisimple
or reductive. We then propose to consider the class-$\gamma$ nilpotent Lie groups of a given dimension all
together, by using the variety of nilpotent Lie algebras, as in the study of compatible metrics (see Section
\ref{varint}). The advantage of this unified approach is that the group giving the isomorphism is the reductive
Lie group $G_{\gamma}$; the price to pay is that the space $\nca_{\gamma}$ where $G_{\gamma}$ is acting on, is
really wild. Fortunately, $\nca_{\gamma}$ is at least a real algebraic variety, and so the classification
problem for such structures may be approached by using tools from invariant theory (see Section \ref{qgv} for an
overview).

We may view $\nca_{\gamma}$ as the variety of all class-$\gamma$ $n$-dimensional nilpotent Lie groups by
identifying each element $\mu\in\nca_{\gamma}$ with a class-$\gamma$ nilpotent Lie group,
\begin{equation}\label{ide1}
\mu\longleftrightarrow   (N_{\mu},\gamma).
\end{equation}
Two class-$\gamma$ structures $\mu,\lambda$ are isomorphic if and only if they live in the same
$G_{\gamma}$-orbit, and hence the quotient
$$
\nca_{\gamma}/G_{\gamma}
$$
parameterizes the moduli space of all $n$-dimensional class-$\gamma$ nilpotent Lie groups up to isomorphism.
Recall that a nilpotent Lie group $N_{\mu}$ admits a class-$\gamma$ structure if and only if the orbit
$\Gl(n).\mu$ meets the variety $\nca_{\gamma}$.

An orbit $G_{\gamma}.\mu$, $\mu\ne 0$, can never be closed since this is equivalent to $(N_{\mu},\gamma)$
admitting a compatible metric satisfying the strong property (\ref{unconditionint}), and thus the categorical
quotient $\nca_{\gamma}\mum G_{\gamma}$ consists only of $(\RR^n,\gamma)$.  Due to the absence of closed orbits,
which are precisely the zeroes of the moment map, it is natural to consider the wider quotient
$\nca_{\gamma}\kir G_{\gamma}$ parameterizing orbits containing a critical point of the functional square norm
of the moment map. Recall that the moment map $m_{\gamma}:\nca_{\gamma}\mapsto\pg_{\gamma}$ for the action of
$G_{\gamma}$ on $\nca_{\gamma}$ is given by $m_{\gamma}(\mu)=8\Ricg_{\mu}$, and therefore $\nca_{\gamma}\kir
G_{\gamma}$ classifies precisely those class-$\gamma$ nilpotent Lie groups admitting a minimal compatible
metric.  Since the negative gradient flow of $F=||\riccig||^2$ stays in the $G_{\gamma}$-orbit of the starting
point, every class-$\gamma$ nilpotent Lie group degenerates via such a flow into one of these special
structures.

We now describe the decomposition of $\nca_{\gamma}\kir G_{\gamma}$ into finitely many categorical quotients due
to L. Ness \cite{Nss}.  Let $[\mu]\in \PP \nca_{\gamma}$ be a critical point of $F$, with
$\Ricg_{\mu}=c_{\mu}I+D_{\mu}$ for some $c_{\mu}\in\RR$ and $D_{\mu}\in\Der(\mu)$.  Then there exists $c>0$ such
that the eigenvalues of $cD_{\mu}$ are all integers prime to each other, say $k_1<...<k_r\in\ZZ$ with
multiplicities $d_1,...,d_r\in\NN$.  The data set $\eigen$ is called the {\it type} of the critical point
$[\mu]$.  The set of types of critical points is in bijection with the finite set of strata determined by the
negative gradient flow of $F$, and it will be denoted by $\{\alpha_1,...,\alpha_s\}$.  For a fixed type $
\alpha=\alpha_i=\eigen$, consider
$$
V_{\alpha}:=\{\mu\in V:D_{\alpha}\in\Der(\mu)\}, \qquad D_{\alpha}=\left[\begin{smallmatrix}
k_1I_{d_1}&&\\
&\ddots&\\
&&k_rI_{d_r}
\end{smallmatrix}\right],
$$
and the reductive Lie group given by
$$
\tilde{G}_{\alpha}:=\left\{ g\in
G_{\gamma}\cap(\Gl(d_1)\times...\times\Gl(d_r)):\prod_{i=1}^{r}(\det{g_i})^{k_i}=1=\det{g}\right\}.
$$
The quotient $\nca_{\gamma}\kir G_{\gamma}$ decomposes as a disjoint union of semi-algebraic varieties
$$
\nca_{\gamma}\kir G_{\gamma}=X_1\cup...\cup X_s,
$$
where each $X_i$ is homeomorphic to the categorical quotient $(V_{\alpha_i}\cap\nca_{\gamma})\mum
\tilde{G}_{\alpha_i}$.  This allows us to approach the classification of invariant geometric structures on
nilpotent Lie groups by using invariant-theoretic methods.  By considering each $X_i$ separately, we have for
instance that geometric structures of type $\alpha_i$ are precisely the closed $\tilde{G}_{\alpha_i}$-orbits,
and two different orbits give rise to non-isomorphic structures.  We therefore have that two non-isomorphic
structures can always been separated by a $\tilde{G}_{\alpha_i}$-invariant polynomial on $V_{\alpha_i}$.
Moreover, $X_i$ can be described by using a set of generators and relations of
$\RR[V_{\alpha_i}\cap\nca_{\gamma}]^{\tilde{G}_{\alpha_i}}$, the ring of all invariant polynomials.  It is shown
that some of the simplest types already lead to wide open problems in invariant theory.

We do not know how far is $\nca_{\gamma}\kir G_{\gamma}$ from the whole quotient $\nca_{\gamma}/G_{\gamma}$. The
crucial question is how strong is, for a class-$\gamma$ structure, the property of admitting a minimal
compatible metric (see end of Section \ref{varint}).

\section{Geometric structures and compatible metrics}\label{geometric}

Let $N$ be a real $n$-dimensional nilpotent Lie group with Lie algebra $\ngo$, whose Lie bracket is denoted by
$\mu :\ngo\times\ngo\mapsto\ngo$.  An invariant geometric structure on $N$ is defined by left translation of a
tensor $\gamma$ on $\ngo$ (or a set of tensors), usually non-degenerate in some way, which satisfies a suitable
integrability condition
\begin{equation}\label{closedg}
{\rm IC}(\gamma,\mu)=0,
\end{equation}
involving only $\mu$ and $\gamma$.  In this paper, we will focus on the following classes of geometric
structures: symplectic, complex and hypercomplex, as well as on their respective `almost' versions, that is,
when condition (\ref{closedg}) is not required.  In this way, ${\rm IC}(\gamma,\mu)$ can be for instance the
differential of a $2$-form or the Nijenhuis tensor associated to some $(1,1)$-tensor. The contact case is
somewhat different because the condition is `open', but it becomes an equation of the form (\ref{closedg}) when
one considers fixed the underlying almost-contact structure.  We shall deal with contact and complex symplectic
structures in a forthcoming paper.

The pair $(N,\gamma)$ will often be called a {\it class-$\gamma$ nilpotent Lie group}, and $N$ will be assumed
to be simply connected for simplicity.  The group $\Gl(n):=\Gl(n,\RR)=\Gl(\ngo)$ of invertible maps of $\ngo$ acts on
the vector space of tensors on $\ngo$ of a given class, preserving the non-degeneracy, and if $\gamma$ is
integrable then $\vp.\gamma$ is so for any $\vp\in\Aut(\ngo)$, the group of automorphisms of $\ngo$. In view of
this fact, two class-$\gamma$ nilpotent Lie groups $(N,\gamma)$ and $(N',\gamma')$ are said to be {\it
isomorphic} if there exists a Lie algebra isomorphism $\vp:\ngo\mapsto\ngo'$ such that $\gamma'=\vp.\gamma$.
Also, given two geometric structures $\gamma,\gamma'$ of the same class on $N$, we say that $\gamma$ {\it
degenerates to} $\gamma'$ if $\gamma'\in\overline{\Aut(\ngo).\gamma}$, the closure of the orbit
$\Aut(\ngo).\gamma$ relative to the natural topology.

A left invariant Riemannian metric on $N$ is said to be {\it compatible} with $(N,\gamma)$ if the corresponding
inner product $\ip$ on $\ngo$ satisfies an orthogonality condition
\begin{equation}\label{ortcon}
{\rm OC}(\gamma,\ip)=0,
\end{equation}
in which only $\ip$ and $\gamma$ are involved.  We denote by $\cca=\cca(N,\gamma)$ the set of all left invariant
metrics on $N$ which are compatible with $(N,\gamma)$.  The pair $(\gamma,\ip)$ with $\ip\in\cca$ will often be
referred to as a {\it class-$\gamma$ metric structure}.  It is clear from (\ref{ortcon}) that for an invariant
geometric structure there always exist a compatible metric, since the condition is independent from $\mu$.
Moreover, the space $\cca$ is usually huge; recall for instance that the group
$$
G_{\gamma}=\{ \vp\in\Gl(n):\vp.\gamma=\gamma\}
$$
acts on $\cca$, and it is easy to see that actually for any $\ip\in\cca$ we have that
\begin{equation}\label{compg}
\cca=G_{\gamma}.\ip=\{\la\vp^{-1}\cdot,\vp^{-1}\cdot\ra:\vp\in G_{\gamma}\}.
\end{equation}
A natural question takes place:
\begin{itemize}
\item[] Given a class-$\gamma$ nilpotent Lie group $(N,\gamma)$, are there distinguished left invariant
Riemannian metrics on $N$ compatible with $\gamma$?
\end{itemize}
As usual, the meaning of the word `distinguished' is part of the question.  This problem may be (and it is) stated
for differentiable manifolds in general, and does not only present some interest in Riemannian geometry; indeed,
the existence of a certain nice compatible metric could eventually help to distinguish two geometric structures
as well as to find privileged geometric structures on a given manifold.

The aim of this section is to propose two properties which make a compatible metric very distinguished, one is
obtained by minimizing a curvature functional and the other as a soliton solution for a natural evolution flow.
The Ricci tensor will be used in both approaches.  In Appendix \ref{app}, we have reviewed some well known
properties of left invariant metrics on nilpotent Lie groups, which will be used constantly from now on.

Fix a class-$\gamma$ nilpotent Lie group $(N,\gamma)$.  Let $\ggo_{\gamma}$ be the Lie algebra of $G_{\gamma}$,
$$
\ggo_{\gamma}=\{ A\in\glg(n):A.\gamma=0\}.
$$

\begin{definition}\label{invricci}
{\rm For each compatible metric, we consider the orthogonal projection $\Ricg_{\ip}$ of the Ricci operator
$\Ric_{\ip}$ on $\ggo_{\gamma}$,  called the {\it invariant Ricci operator}, and the corresponding {\it
invariant Ricci tensor} given by $\riccig=\la\Ricg\cdot,\cdot\ra$. }
\end{definition}

The role of the Ricci tensor has always been crucial in defining privileged (compatible) metrics; we have for
example Einstein metrics, extremal K$\ddot{{\rm a}}$hler metrics in complex geometry, and more recently metrics
with hermitian Ricci tensor and $\phi$-invariant Ricci tensor in symplectic and contact geometry, respectively.
These two last notions are equivalent to $\riccig=0$ and have been characterized in the compact case by D.
Blair, S. Ianus and A. Ledger \cite{BlrIns, BlrLdg} as the critical points of two very natural
curvature functionals on $\cca$: the total scalar curvature functional $S$ and a functional $K$ for which the
global minima are precisely K$\ddot{{\rm a}}$hler or Sasakian metrics, respectively (see also \cite{Blr} and
Section \ref{hermricci}).

In this light, condition
\begin{equation}\label{uncondition}
\riccig_{\ip}=0,
\end{equation}
involves both the geometric structure and the metric, and seems to be very natural to require to a compatible metric.
Nevertheless, if $\RR I\subset\ggo_{\gamma}$, then $\tr{\Ricg_{\ip}}=\scalar(\ip)$, and so it is forbidden for
instance for non-abelian nilpotent Lie groups (where always $\scalar(\ip)<0$) in the complex and hypercomplex
cases. We shall prove that this condition is forbidden in the symplectic case as well.  We therefore have to
consider (\ref{uncondition}) as an unreachable goal and try to get as close as possible, for instance, by
minimizing $||\riccig_{\ip}||^2=\tr(\Ricg_{\ip})^2$. In order to avoid homothetical changes, we must normalize
the metrics some way.  In the noncompact homogeneous case, the scalar curvature always appears as a very natural
choice.  We then propose the following

\begin{definition}\label{minimalg}
{\rm A left invariant metric $\ip$ compatible with a class-$\gamma$ nilpotent Lie group $(N,\gamma)$ is called
{\it minimal} if
$$
\tr(\Ricg_{\ip})^2=\min\{  \tr(\Ricg_{\ip'})^2 : \ip'\in\cca(N,\gamma), \quad \scalar(\ip')=\scalar(\ip)\}.
$$ }
\end{definition}

Recall that the existence and uniqueness (up to isometry and scaling) of minimal metrics is far to be clear from
the definition.  The uniqueness shall be proved in Section \ref{var}, but the `existence question' is still
nebulous.  Minimal metrics are the compatible metrics closest to satisfy the `Einstein-like' condition
$\riccig_{\ip}=c\ip$, for some $c\in\RR$.   Indeed,
$$
||\Ricg_{\ip}-\frac{\tr{\Ricg_{\ip}}}{n}I||^2=\tr(\Ricg_{\ip})^2-\frac{(\tr{\Ricg_{\ip}})^2}{n}
$$
and $\tr{\Ricg_{\ip}}$ equals either $0$ or $\scalar(\ip)$, depending on $\ggo_{\gamma}$ contains or not $\RR I$.

We now consider an evolution approach.  Motivated by the famous Ricci flow introduced by R. Hamilton
\cite{Hml1}, we consider the {\it invariant Ricci flow} for our left invariant metrics on $N$, given by the
following evolution equation
\begin{equation}\label{grfn}
\ddt \ip_t=\pm\riccig_{\ip_t},
\end{equation}
which coincides for example with the anticomplexified  Ricci flow studied in \cite{LeWng} in the symplectic case
(see Section \ref{hermricci}).  The choice of the best sign might depend on the class of structure.  This is
just an ordinary differential equation and hence the existence for all $t$ and uniqueness of the solution is
guaranteed. It follows from (\ref{compg}) that
\begin{equation}\label{tangcompg}
\tang_{\ip}\cca=\{\alpha\in\sym(\ngo):A_{\alpha}.\gamma=0\},
\end{equation}
and therefore, if $\ip_0\in\cca$ then the solution $\ip_t\in\cca$ for all $t$ since
$\riccig_{\ip_t}\in\tang_{\ip_t}\cca$ (see Appendix \ref{app}).

\begin{remark}\label{uniq}
{\rm If we had however the uniqueness of the solution for the flow (\ref{grfn}) in the non-compact general case,
then we would not need to restrict ourselves to left invariant metrics. Indeed, if $f$ is an isometry of the
initial metric $\ip_0$ which also leaves $\gamma$ invariant, then since $f^*\ip_t$ is also a solution and
$f^*\ip_0=\ip_0$ we would get by uniqueness of the solution that $f$ is an isometry of all the metrics $\ip_t$
as well. Left invariance of the starting metric would be therefore preserved along the flow. }
\end{remark}

When $M$ is compact, a normalized Ricci flow is often considered, under which the volume of the solution metric
is constant in time.  Actually, the normalized equation differs from the Ricci flow only by a change of scale in
space and a change of parametrization in time (see \cite{Hml2, CaoChw}).  In our case, where the manifold is
non-compact but the metrics are homogeneous, it seems natural to do the same thing but normalizing by the scalar
curvature, which is just a single number associated to the metric.  We recall that a left invariant metric $\ip$
on a nilpotent Lie group $N$ has always $\scalar(\ip)<0$, unless $N$ is abelian (see (\ref{ricci})).

\begin{proposition}\label{nacrf}
The solution to the normalized invariant Ricci flow
\begin{equation}\label{nfg}
\ddt\ip_t=\pm\riccig_{\ip_t}\mp\frac{\tr(\Ricg_{\ip_t})^2}{\scalar(\ip_t)}\ip_t
\end{equation}
satisfies $\scalar(\ip_t)=\scalar(\ip_0)$ for all $t$.  Moreover, this flow differs from the invariant Ricci
flow {\rm (\ref{grfn})} only by a change of scale in space and a change of parametrization in time.
\end{proposition}

\begin{proof}
It follows from the formula for the gradient of the scalar curvature functional $\scalar:\pca\mapsto\RR$ given
in (\ref{gradsc}) that if $\ip_t$ is a solution of (\ref{nfg}), then the function $f(t)=\scalar(\ip_t)$
satisfies
$$
\begin{array}{rl}
f'(t)&=g_{\ip_t}(\ddt\ip_t,-\ricci_{\ip_t})\\ \\
&=\mp g_{\ip_t}(\riccig_{\ip_t},\ricci_{\ip_t})\pm\frac{\tr(\Ricg_{\ip_t})^2}{\scalar(\ip_t)}
g_{\ip_t}(\ip_t,\ricci_{\ip_t})\\ \\
&=\mp\tr(\Ricg_{\ip_t}\Ric_{\ip_t})\pm\frac{\tr(\Ricg_{\ip_t})^2}{\scalar(\ip_t)}\tr(\Ric_{\ip_t})\\ \\
&=\mp\tr(\Ricg_{\ip_t})^2(1-\frac{f(t)}{\scalar(\ip_t)})=0, \qquad \forall t,
\end{array}
$$
and thus $f(t)\equiv f(0)=\scalar(\ip_0)$.  The last assertion follows as in \cite{Hml2} in a completely
analogous way.
\end{proof}

The fixed points of this normalized flow (\ref{nfg}) are those metrics satisfying $\Ricg_{\ip}\in\RR I$, and so
in particular, if $\ggo_{\gamma}\subset\slg(n)$, then this is equivalent to $\Ricg_{\ip}=0$.  Indeed, if
$\Ricg_{\ip}=\mp\frac{\tr(\Ricg_{\ip_t})^2}{\scalar(\ip_0)}I$ then $\Ricg_{\ip}=0$ since $\tr\Ricg_{\ip}=0$.  We
should also note that for the flow (\ref{nfg}), $\ddt\ip_t\in\tang_{\ip_t}\cca+\RR I$ for all $t$, which implies
that the solution $\ip_t$ stays in the set of all scalar multiples of compatible metrics.  Recall that if $\RR
I\subset\ggo_{\gamma}$ then the solution stays anyway in $\cca$.

In these evolution approaches always appear naturally the soliton metrics, which are not fixed points of the
flow but are close to, and they play an important role in the study of singularities (see the surveys
\cite{Hml2, CaoChw} for further information).  The idea is that if one is trying to improve a metric via an
evolution equation, then those metrics for which the solution remains isometric to the initial point may be
certainly considered as very distinguished.

\begin{definition}\label{solitong}
{\rm A metric $\ip$ compatible with $(N,\gamma)$ is called an {\it invariant Ricci soliton} if the solution
$\ip_t$ to the normalized invariant Ricci flow (\ref{nfg}) with initial metric $\ip_0=\ip$ is given by
$\vp_t^*\ip$, the pullback of $\ip$ by a one parameter group of diffeomorphisms $\{\vp_t\}$ of $N$. }
\end{definition}

We now give a neat characterization of invariant Ricci soliton metrics, which will be very useful in Section
\ref{var} to prove the equivalence with the property of being minimal (see Definition \ref{minimalg}), and to
find explicit examples in the subsequent sections.

\begin{proposition}\label{soliton1g}
Let $(N,\gamma)$ be a class-$\gamma$ nilpotent Lie group.  A compatible metric $\ip$ is an invariant Ricci
soliton if and only if $\Ricg_{\ip}=cI+D$ for some $c\in\RR$ and $D\in\Der(\ngo)$.  In such a case,
$c=\frac{\tr(\Ricg_{\ip})^2}{\scalar(\ip)}$.
\end{proposition}

\begin{proof}
We first note that the assertion on the value of the number $c$ follows from (\ref{ricort}); in fact,
$$
\tr(\Ricg_{\ip})^2=\tr(\Ric_{\ip}\Ricg_{\ip})=c\tr{\Ric_{\ip}}+\tr(\Ric_{\ip}D)=c\scalar(\ip).
$$
Assume that there exists a one-parameter group of diffeomorphisms $\vp_t$ on $N$ such that $\ip_t=\vp^*_t\ip$ is
a solution to the flow (\ref{nfg}).  By the uniqueness of the solution we have that $\vp_t^*g$ is also
$N$-invariant for all $t$ (see Remark \ref{uniq}).  Thus $\vp_t$ normalizes $N$ and so it follows from \cite[Thm
2, 4)]{Wls} that $\vp_t\in\Aut(N).N$.  This implies that there exists a one-parameter group $\psi_t$ of
automorphisms of $N$ such that $\vp_t^*\ip=\psi_t^*\ip$ for all $t$.  Now, if $\psi_t=e^{-\frac{t}{2}D}$ with
$D\in\Der(\ngo)$ then $\ddt|_0\psi_t^*\ip=\la D\cdot,\cdot\ra$, and using that $\psi_t^*\ip$ is a solution of
(\ref{nfg}) in $t=0$ we obtain that $\riccig_{\ip}=c\ip+\la D\cdot,\cdot\ra$ for some $c\in\RR$, or
equivalently, $\Ricg=cI+D$.

Conversely, if $\Ricg=cI+D$ then we will show that the curve $\ip_t=e^{-\frac{t}{2}D}.\ip$ is a solution of the
flow (\ref{nfg}).  For any $t$, it follows from
$\frac{t}{2}D=\frac{t}{2}\Ricg_{\ip}-\frac{t}{2}cI\in\ggo_{\gamma}+\RR I$ that
$$
\gamma=b(t)e^{-\frac{t}{2}D}.\gamma.
$$
for some $b(t)\in\RR$.  This implies that
$$
\Ricg_{\ip_t}=e^{-\frac{t}{2}D}\Ricg_{\ip}e^{\frac{t}{2}D}=e^{-\frac{t}{2}D}(cI+D)e^{\frac{t}{2}D}=cI+D
$$
for all $t$.  Therefore
$$
\begin{array}{rl}
\ddt|_0\ip_t&=\la D\cdot,\cdot\ra_t=\la(\Ricg_{\ip_t}-cI)\cdot,\cdot\ra_t \\ \\
&=\riccig_{\ip_t}-c\ip_t=\riccig_{\ip_t}+\frac{\tr(\Ricg_{\ip_t})^2}{\scalar(\ip_0)}\ip_t,
\end{array}
$$
as was to be shown.
\end{proof}

Recall that the condition in the above proposition can be replaced by
$$
\Ricg_{\ip}-\frac{\tr(\Ricg_{\ip})^2}{\scalar(\ip)}I\in\Der(\ngo),
$$
which gives a computable method to check whether a metric is an invariant Ricci soliton or not.

\section{Real geometric invariant theory and the moment map}\label{git}

In this section, we overview some results from (geometric) invariant theory over the real numbers.  We refer to
\cite{RchSld} for a detailed exposition.  These will be our tools to study metrics compatible with geometric
structures on nilpotent Lie groups, as well as to approach the classification problem for such structures.

\subsection{Closed orbits and minimal vectors}\label{comv}

Let $G$ be a real reductive Lie group acting on a real vector space $V$ (see \cite{RchSld} for a precise
definition of the situation).  Let $\ggo$ denote the Lie algebra of $G$ with Cartan decomposition
$\ggo=\kg\oplus\pg$, where $\kg$ is the Lie algebra of a maximal compact subgroup $K$ of $G$.  Endow $V$ with a
fixed from now on inner product $\ip$ such that $\kg$ and $\pg$ act by skew-symmetric and symmetric
transformations, respectively. Let $\mca=\mca(G,V)$ denote the set of {\it minimal vectors}, that is
$$
\mca=\{ v\in V: ||v||\leq ||g.v||\quad\forall g\in G\}.
$$
For each $v\in V$ define
$$
\rho_v:G\mapsto\RR, \qquad \rho_v(g)=||g.v||^2=\la g.v,g.v\ra.
$$
In \cite{RchSld}, R. Richardson and P. Slodowy showed that the nice interplay between closed orbits and minimal
vectors found by G. Kempf and L. Ness for actions of complex reductive algebraic groups, is still valid in the
real situation.

\begin{theorem}\cite{RchSld}\label{RS}
Let $V$ be a real representation of a real reductive Lie group $G$, and let $v\in V$.
\begin{itemize}
\item[(i)] $v\in\mca$ if and only if $\rho_v$ has a critical point at $e\in G$.

\item[(ii)] If $v\in\mca$ then $G.v\cap\mca=K.v$.

\item[(iii)] The orbit $G.v$ is closed if and only if $G.v$ meets $\mca$.

\item[(iv)] The closure of any orbit $G.v$ always meets $\mca$.  Indeed, there exists $A\in\pg$ such that
$\lim_{t\to-\infty}\exp(tA).v=v_0$ exists and $G.v_0$ is closed.

\item[(v)] $\overline{G.v}\cap\mca$ is a single $K$-orbit, or in other words, $\overline{G.v}$ contains a unique
closed $G$-orbit.
\end{itemize}
\end{theorem}

As usual in the real case, classical topology of $V$ is always considered rather than Zarisky topology, unless
explicitly indicated otherwise.

\subsection{Critical points of the moment map}\label{cpmp}

We keep the notation of the above subsection.  Let $(\dif\rho_v)_e:\ggo\mapsto\RR$ denote the differential of
$\rho_v$ at the identity $e$ of $G$. It follows from the $K$-invariance of $\ip$ that $(\dif\rho_v)_e$ vanishes
on $\kg$, and so we can assume that $(\dif\rho_v)_e\in\pg^*$, the vector space of real-valued functionals on
$\pg$.  We therefore may define a function called the {\it moment map} for the action of $G$ on $V$ by
\begin{equation}\label{moment}
m:V\mapsto\pg, \qquad  \la m(v),A\ra_{\pg}=(\dif\rho_v)_e(A),
\end{equation}
where $\ip_{\pg}$ is an $\Ad(K)$-invariant inner product on $\pg$.  Since $m(tv)=t^2m(v)$ for all $t\in\RR$, we
also may consider the moment map
\begin{equation}\label{moment2}
m:\PP V\mapsto\pg, \qquad  m(x)=\frac{m(v)}{||v||^2}, \quad 0\ne v\in V,\;x=[v],
\end{equation}
where $\PP V$ is the projective space of lines in $V$.  If $\pi:V\setminus\{ 0\}\mapsto \PP V$ denotes the usual
projection map, then $\pi(v)=x$.  In the complex case, under the natural identifications
$\pg=\pg^*=(\im\kg)^*=\kg^*$, the function $m$ is precisely the moment map from symplectic geometry,
corresponding to the Hamiltonian action of $K$ on the symplectic manifold $\PP V$ (see for instance the survey
\cite{Krw2} or \cite[Chapter 8]{Mmf} for further information).

Consider the functional square norm of the moment map
\begin{equation}\label{norm}
F:V\mapsto\RR, \qquad  F(v)=||m(v)||^2=\la m(v),m(v)\ra_{\pg},
\end{equation}
which is easily seen to be a $4$-degree homogeneous polynomial.  Recall that $\mca$ coincides with the set of
zeros of $F$. It then follows from Theorem \ref{RS}, parts (i) and (iii), that an orbit $G.v$ is closed if and
only if $F(w)=0$ for some $w\in G.v$, and in that case, the set of zeros of $F|_{G.v}$ coincides with $K.v$.  We
furthermore have the following result due to A. Neeman and G. Schwarz (see \cite{RchSld}).

\begin{theorem}\label{retraction}
Let $X$ be a closed $G$-invariant subset of $V$ and set $\mca_X=\mca\cap X$.  Then the negative gradient flow of
$F:V\mapsto\RR$ defines a $K$-equivariant deformation retraction $\psi:X\times[0,1]\mapsto X$ from $X$ onto
$\mca_X$ along $G$-orbits, that is,
\begin{itemize}
\item[(i)] $\psi_0=Id$, $\psi_1(X)=\mca_X$, $\psi(v,t)=v$ for any $v\in\mca_X$.

\item[(ii)] $\psi(v,t)\in G.v$ for all $t<1$ and so $\psi(v,1)\in\overline{G.v}$ for any $v\in X$.
\end{itemize}
$\psi$ also determines a deformation retraction of $X/G$ onto $\mca_X/K$.
\end{theorem}

A natural question arises: what is the role played by the remaining critical points of $F:\PP V\mapsto\RR$ (i.e.
those for which $F(x)>0$) in the study of the $G$-orbit space of the action of $G$ on $V$, as well as on other
real $G$-varieties contained in $V$?. This was studied independently by F. Kirwan \cite{Krw1} and L. Ness
\cite{Nss}, and it is shown in the complex case that the non-minimal critical points share some of the nice
properties of minimal vectors stated in Theorem \ref{RS}.  In the real case, which is actually the one we need
to apply in this paper, the analogous of some of these results have been proved by A. Marian \cite{Mrn}.

 \begin{theorem}\cite{Mrn}\label{marian}
Let $V$ be a real representation of a real reductive Lie group $G$,  $m:\PP V\mapsto\pg$ the moment map and
$F=||m||^2:\PP V\mapsto\RR$.
\begin{itemize}
\item[(i)] If $x\in\PP V$ is a critical point of $F$ then the functional $F|_{G.x}$ attains its minimum value at
$x$.

\item[(ii)] If nonempty, the critical set of $F|_{G.x}$ consists of a unique $K$-orbit.
\end{itemize}
\end{theorem}

We endow $\PP V$ with the Fubini-Study metric defined by $\ip$ and denote by $x\mapsto A_x$ the vector field on
$\PP V$ defined by $A\in\ggo$ via the action of $G$ on $\PP V$, that is, $A_x=\ddt|_0\exp(tA).x$.

\begin{lemma}\cite{Mrn}\label{marian2}
The gradient of the functional $F=\parallel m\parallel^2:\PP V\mapsto\RR$ is given by
$$
\grad(F)_x=4m(x)_x, \qquad x\in\PP V,
$$
and therefore $x$ is a critical point of $F$ if and only if $m(x)_x=0$, if and only if $\exp{tm(x)}$ fixes $x$.
\end{lemma}

We now develop some examples which are far to be the natural ones, but they are those ones will be considered in
this paper to study left invariant structures on nilpotent Lie groups.

\begin{example}\label{gln}
{\rm Let $\ngo$ be an $n$-dimensional real vector space and $V=\lam$ the vector space of all skew-symmetric
bilinear maps from $\ngo\times\ngo$ to $\ngo$. There is a natural action of $\Gl(n):=\Gl(n,\RR)$ on $V$ given by
\begin{equation}\label{action}
g.\mu(X,Y)=g\mu(g^{-1}X,g^{-1}Y), \qquad X,Y\in\ngo, \quad g\in\Gl(n),\quad \mu\in V.
\end{equation}
Any inner product $\ip$ on $\ngo$ defines an $\Or(n)$-invariant inner product on $V$, denoted also by $\ip$, as
follows:
\begin{equation}\label{innV}
\la\mu,\lambda\ra=\sum_{ijk}\la\mu(X_i,X_j),X_k\ra\la\lambda(X_i,X_j),X_k\ra,
\end{equation}
where $\{ X_1,...,X_n\}$ is any orthonormal basis of $\ngo$.  A Cartan decomposition for the Lie algebra of
$\Gl(n)$ is given by $\glg(n)=\sog(n)\oplus\sym(n)$, that is, in skew-symmetric and symmetric transformations
respectively, and we consider the following $\Ad(\Or(n))$-invariant inner product on $\pg:=\sym(n)$,
\begin{equation}\label{inng}
\la A,B\ra_{\pg}=\tr{AB}, \qquad A,B\in\pg.
\end{equation}
The action of $\glg(n)$ on $V$ obtained by differentiation of (\ref{action}) is given by
\begin{equation}\label{actiong}
A.\mu=-\delta_{\mu}(A):= A\mu(\cdot,\cdot)-\mu(A\cdot,\cdot)-\mu(\cdot,A\cdot), \qquad A\in\glg(n),\quad\mu\in
V.
\end{equation}
If $\mu\in V$ satisfies the Jacobi condition, then $\delta_{\mu}:\glg(n)\mapsto V$ coincides with the cohomology
coboundary operator of the Lie algebra $(\ngo,\mu)$ from level $1$ to $2$, relative to cohomology with values in
the adjoint representation.  Recall that $\Ker{\delta_{\mu}}=\Der(\mu)$, the Lie algebra of derivations of the
algebra $\mu$.  Let $A^t$ denote the transpose relative to $\ip$ of a linear transformation $A:\ngo\mapsto\ngo$
and consider the adjoint map $\ad_{\mu}{X}:\ngo\mapsto\ngo$ (or left multiplication) defined by
$\ad_{\mu}{X}(Y)=\mu(X,Y)$.

\begin{proposition}\label{ourmoment}
The moment map $m:V\mapsto\pg$ for the action {\rm (\ref{action})} of $\Gl(n)$ on $V=\lam$ is given by
\begin{equation}\label{defm}
m(\mu)=-4\displaystyle{\sum_{i}}(\ad_{\mu}{X_i})^t\ad_{\mu}{X_i}
+2\displaystyle{\sum_{i}}\ad_{\mu}{X_i}(\ad_{\mu}{X_i})^t,
\end{equation}
where $\{ X_1,...,X_n\}$ is any orthonormal basis of $\ngo$, and it is a simple calculation to see that
\begin{equation}\label{defm2}
\begin{array}{rl}
\la m(\mu)X,Y\ra=&-4\displaystyle{\sum_{ij}}\la\mu(X,X_i),X_j\ra\la\mu(Y,X_i),X_j\ra \\
&+2\displaystyle{\sum_{ij}}\la\mu(X_i,X_j),X\ra\la\mu(X_i,X_j), Y\ra, \qquad \forall\; X,Y\in\ngo.
\end{array}
\end{equation}
\end{proposition}

\begin{proof}
For any $A\in\pg$ we have that
$$
\begin{array}{rl}
(\dif\rho_{\mu})_I(A)=&\ddt|_0\la e^{tA}.\mu,e^{tA}.\mu\ra =-2\la\mu,\delta_{\mu}(A)\ra \\ \\

=&-2\displaystyle{\sum_{pij}}
\la\mu(X_p,X_i),X_j\ra\la\delta_{\mu}(A)(X_p,X_i),X_j\ra \\

=&-2\Big(\displaystyle{\sum_{pij}}\la\mu(X_p,X_i),X_j\ra\la\mu(AX_p,X_i),X_j\ra \\

&+ \la\mu(X_p,X_i),X_j\ra\la\mu(X_p,AX_i),X_j\ra  \\ \\

&-\la\mu(X_p,X_i),X_j\ra\la A\mu(X_p,X_i),X_j\ra\Big) \\ \\

=&
-2\Big(\displaystyle{\sum_{pijr}}\la\mu(X_p,X_i),X_j\ra\la\mu(X_r,X_i),X_j\ra\la AX_p,X_r\ra \\

& +\la\mu(X_p,X_i),X_j\ra\la\mu(X_p,X_r),X_j\ra\la AX_i,X_r\ra  \\ \\

& -\la\mu(X_p,X_i),X_j\ra\la\mu(X_p,X_i),X_r\ra\la AX_j,X_r\ra\Big).
\end{array}
$$
By interchanging the indexes $p$ and $i$ in the second line, and $p$ and $j$ in the third one, we get
$$
\begin{array}{rl}
(\dif\rho_{\mu})_I(A)=& -4\displaystyle{\sum_{prij}}\la\mu(X_p,X_i),X_j\ra \la\mu(X_r,X_i),X_j\ra
\la AX_p,X_r\ra\\

&+2\displaystyle{\sum_{prij}}\la\mu(X_i,X_j),X_p\ra\la\mu(X_i,X_j),X_r\ra\la AX_p,X_r\ra.
\end{array}
$$
If we call $M$ the right hand side of (\ref{defm}), then we obtain from (\ref{defm2}) that
$$
(\dif\rho_{\mu})_I(A)=\displaystyle{\sum_{pr}} \la MX_p,X_r\ra\la AX_p,X_r\ra=\tr{MA}=\la M,A\ra_{\pg},
$$
which implies that $m(\mu)=M$ (see (\ref{moment})).
\end{proof}   }
\end{example}

\begin{example}\label{momg}
{\rm We keep the notation as in Example \ref{gln}.  Let $\gamma$ be a tensor on $\ngo$ and let
$G_{\gamma}\subset\Gl(n)$ denote the subgroup leaving $\gamma$ invariant, with Lie algebra $\ggo_{\gamma}$.  The
group $G_{\gamma}$ is reductive and $K_{\gamma}=G_{\gamma}\cap\Or(n)$ is a maximal compact subgroup of
$G_{\gamma}$, whose Lie algebra will be denoted by $\kg_{\gamma}$.  A Cartan decomposition is given by
$$
\ggo_{\gamma}=\kg_{\gamma}\oplus\pg_{\gamma}, \qquad \pg_{\gamma}=\pg\cap\ggo_{\gamma}.
$$
If $p:\pg\mapsto\pg_{\gamma}$ is the orthogonal projection relative to $\ip_{\pg}$, then it is easy to see that
$$
m_{\gamma}:V\mapsto\pg_{\gamma}, \qquad m_{\gamma}=p\circ m,
$$
is precisely the moment map for the action of $G_{\gamma}$ on $V$.   }
\end{example}

In the cases considered in detail in this paper we will have $(G_{\gamma},K_{\gamma})$ equal to
$(\Spe(\frac{n}{2},\RR),\U(\frac{n}{2}))$ (symplectic), $(\Gl(\frac{n}{2},\CC),\U(\frac{n}{2}))$ (complex),
$(\Gl(\frac{n}{4},\HH),\Spe(\frac{n}{4}))$ (hypercomplex) and $(\Gl(n,\RR),\Or(n))$ ($\gamma=0$).

\subsection{Quotients of $G$-varieties}\label{qgv}

Let $G$ be a real reductive Lie group acting on a finite dimensional real vector space $V$, and let $X\subset V$
be a $G$-variety, that is, a real algebraic variety which is $G$-invariant.  The main problem of geometric
invariant theory is to understand the orbit space of the action of $G$ on $X$, parameterized by the quotient
$X/G$ (see \cite{PpvVnb} for further information).  In the real situation, the main references are
\cite{RchSld}, on which this brief overview is based, and papers due to D. Luna and G. Schwarz cited in
\cite{RchSld}.

The standard quotient topology of $X/G$ can be very ugly, for instance, if $y\in\overline{G.x}$ and $G.x\ne G.y$
then they can not be separated by $G$-invariant open neighborhoods and so $X/G$ is not usually a $T_1$-space. In
order to avoid this problem one may consider a smaller quotient $X\mum G$ parame\-tri\-zing only closed orbits.
For any $x\in X$ there exists a unique closed $G$-orbit in $\overline{G.x}$, and so a natural map $q:X\mapsto
X\mum G$ can be defined.  Consider on $X\mum G$ the quotient topology for the map $q$.  $X\mum G$ is called the
{\it categorical quotient} for the action of $G$ on $X$ since it satisfies the following universal property in
the category of all $T_1$-spaces: for any continuous map $\alpha:X\mapsto Y$ that is constant on $G$-orbits,
there exists a unique continuous map $\beta:X\mum G\mapsto Y$ such that $\alpha=\beta\circ q$.  The uniqueness
of an object with such a property is clear from the definition.  Recall that the usual quotient $X/G$ would be
the categorical quotient in the category of all topological spaces.

The space $X\mum G$ is actually Hausdorff.  Moreover, the map $\mca_X/K\mapsto X\mum G$ determined by the
inclusion $\mca_X\mapsto X$ is a homeomorphism, and hence it follows from properties of actions of compact Lie
groups that $X\mum G$ is a real semi-algebraic variety (i.e. defined by finitely many polynomial equations and
inequalities).

It is clear that $X\mum G$ is the `correct' quotient to consider from many points of view, but the price to pay
for a $T_1$ quotient is that in some cases $X\mum G$ classifies only a very few orbits.  In that case, it is
then natural to try to understand the bigger set $X\kir G$ parameterizing those orbits which contain a
(non-necessarily zero) critical point of $F=||m||^2:\PP V\mapsto\RR$.  Let $C\subset\PP V$ denote the set of all
such critical points and set $C_X:=C\cap\PP X$, where $\PP X=\pi(X)$ is the projection of $X$ on $\PP V$.  Since
there is at most one $K$-orbit of critical points on each $G.[v]$, we have that $C_X/K=GC_X/G$ and thus we may
define
$$
X\kir G:=\{ v\in X:[v]\in C_X\}/K=\pi^{-1}(C_X/K), \qquad \PP X\kir G:=C_X/K.
$$
Recall that $X\kir G$ is invariant by scalar multiplication and so the fact that two vectors in a line $\RR^*v$
were in the same $G$-orbit or not is disregarded.

This is not the place for deep considerations on the topology of this wider quotient $X\kir G$ (nor are we
qualified to do it).  We do not know for instance if $X\kir G$ will be the categorical quotient for a suitable
category of non-necessarily $T_1$ topological spaces; or if the quotient topology for the map $q:\PP X\mapsto
\PP X\kir G$ defined by the negative gradient flow of $F$ coincides with the standard topology of $C_X/K$.

Anyway, although the space $X\kir G$ is rather nebulous, the following result due to L. Ness helps us to
describe it, by giving a decomposition of $X\kir G$ in finitely many disjoint subsets each one being a
categorical quotient for a suitable action and so homeomorphic to a semi-algebraic variety.  These results are
proved in the complex case, but it is not hard to see that they remain valid over $\RR$.

\begin{theorem}\cite{Nss}\label{ness2}
The negative gradient flow of $F=||m||^2:\PP V\mapsto\RR$ determines a stratification of the projection $\PP X$
of $X$ on $\PP V$ given by
$$
\PP X=S_{\la A_1\ra}\cup...\cup S_{\la A_s\ra}, \qquad A_1,...,A_s\in\pg,
$$
where each stratum $S_{\la A\ra}$, $A\in\pg$, is the set of all the points $x\in \PP X$ which flow into $C_{\la
A\ra}$, the set of critical points $y\in \PP X$ of $F$ such that $m(y)\in\Ad(K).A$.  Moreover, for each stratum
$S_{\la A\ra}$ there exists a subspace $V_ A\subset V$ and a reductive subgroup $\tilde{G}_A\subset G$ such that
$C_{\la A\ra}/K =\pi((V_A\cap X)\mum \tilde{G}_A$), the projection of the categorical quotient for the action of
$\tilde{G}_{A}$ on $V_{A}\cap X$.
\end{theorem}

It then follows that $X\kir G$ can be decomposed as a disjoint union
\begin{equation}\label{decomp}
X\kir G=\pi^{-1}(C_{\la A_1\ra}/K)\cup...\cup\pi^{-1}(C_{\la A_s\ra}/K), \qquad A_1,...,A_s\in\pg,
\end{equation}
where each subset $\pi^{-1}(C_{\la A_i\ra}/K)$, $i=1,...,s$, is homeomorphic to the categorical quotient
$(V_{A_i}\cap X)\mum \tilde{G}_{A_i}$ and hence is a semi-algebraic variety.

In what follows, we explain how to construct the subspace $V_A$ and the subgroup $\tilde{G}_A$.  For each
$A\in\pg$ consider
$$
V_A=\{ v\in V:A.v={\small \frac{||A||^2}{2}}v\}, \qquad C_A=\{[v]\in C_X:m([v])=A\},
$$
and
$$
G_A=\{ g\in G:\Ad(g)A=A\}, \qquad \ggo_A={\rm Lie}(G_A)=\{ B\in\ggo:[B,A]=0\}.
$$
Let $\tilde{G}_A$ be the subgroup of $G_A$ with Lie algebra $\tilde{\ggo}_A=\ggo_A/\RR A$.  The corresponding
Cartan decompositions will be denoted by $\ggo_A=\kg_A\oplus\pg_A$ and
$\tilde{\ggo}_A=\kg_A\oplus\tilde{\pg}_A$, where $\pg_A=\tilde{\pg}_A\oplus\RR A$ is an orthogonal decomposition
and $\kg_A$ is the Lie algebra of $K_A=\{ g\in K:\Ad(g)A=A\}$, the maximal compact subgroup of both $G_A$ and
$\tilde{G}_A$. A crucial point here proved in \cite{Nss} is that $m(V_A)\subset\pg_A$ and hence the moment map
$m_A:V_A\mapsto\tilde{\pg}_A$ for the action of $\tilde{G}_A$ on $V_A$ is just given by
$$
m_A=p\circ m|_{V_A},
$$
where $p:\pg_A\mapsto\tilde{\pg}_A$ is the orthogonal projection and $m:V\mapsto\pg$ is the moment map for the
action of $G$ on $V$.  This implies that
$$
C_A=\{[v]\in\PP (V_A\cap X):m_A([v])=0\}=\pi(\{ v\in V_A\cap X:m_A(v)=0\}),
$$
and hence
$$
C_{\la A\ra}/K=C_A/K_A=\pi((V_A\cap X)\mum\tilde{G}_A),
$$
the projection of the corresponding categorical quotient.

\section{Variety of compatible metrics}\label{var}

Let us consider as a parameter space for the set of all real nilpotent Lie algebras of a given dimension $n$, the
set $\nca$ of all nilpotent Lie brackets on a fixed $n$-dimensional real vector space $\ngo$.  If
$$
V=\lam=\{\mu:\ngo\times\ngo\mapsto\ngo : \mu\; \mbox{skew-symmetric bilinear map}\},
$$
then
$$
\nca=\{\mu\in V:\mu\;\mbox{satisfies Jacobi and is nilpotent}\}
$$
is an algebraic subset of $V$.  Indeed, the Jacobi identity and the nilpotency condition are both determined by
zeroes of polynomials.

We fix a tensor $\gamma$ on $\ngo$ (or a set of tensors), and let $G_{\gamma}$ denote the subgroup of $\Gl(n)$
preserving $\gamma$. These groups act naturally on $V$ by (\ref{action}) and leave $\nca$ invariant.  Consider
the subset $\nca_{\gamma}\subset\nca$ given by
$$
\nca_{\gamma}=\{\mu\in\nca:{\rm IC}(\gamma,\mu)=0\},
$$
that is, those nilpotent Lie brackets for which $\gamma$ is integrable (see (\ref{closedg})). $\nca_{\gamma}$ is
also an algebraic variety since ${\rm IC}(\gamma,\mu)$ is always linear on $\mu$.  Recall that
$$
W_{\gamma}=\{\mu\in V:{\rm IC}(\gamma,\mu)=0\}
$$
is a $G_{\gamma}$-invariant linear subspace of $V$, and $\nca_{\gamma}=\nca\cap W_{\gamma}$.

For each $\mu\in\nca$, let $N_{\mu}$ denote the simply connected nilpotent Lie group with Lie algebra
$(\ngo,\mu)$.  We now consider an identification of each point of $\nca_{\gamma}$ with a compatible metric.  Fix
an inner product $\ip$ on $\ngo$ compatible with $\gamma$, that is, such that (\ref{ortcon}) holds. We identify
each $\mu\in\nca_{\gamma}$ with a class-$\gamma$ metric structure
\begin{equation}\label{ideg}
\mu\longleftrightarrow   (N_{\mu},\gamma,\ip),
\end{equation}
where all the structures are defined by left invariant translation.  Therefore, each $\mu\in\nca_{\gamma}$ can
be viewed in this way as a metric compatible with the class-$\gamma$ nilpotent Lie group $(N_{\mu},\gamma)$, and
two metrics $\mu,\lambda$ are compatible with the same geometric structure if and only if they live in the same
$G_{\gamma}$-orbit. Indeed, the action of $G_{\gamma}$ on $\nca_{\gamma}$ has the following interpretation: each
$\vp\in G_{\gamma}$ determines a Riemannian isometry preserving the geometric structure
$$
(N_{\vp.\mu},\gamma,\ip)\mapsto (N_{\mu},\gamma,\la\vp\cdot,\vp\cdot\ra)
$$
by exponentiating the Lie algebra isomorphism $\vp^{-1}:(\ngo,\vp.\mu)\mapsto(\ngo,\mu)$.  We then have the
identification $G_{\gamma}.\mu=\cca(N_{\mu},\gamma)$, and more in general the following

\begin{proposition}\label{upto}
Every class-$\gamma$ metric structure $(N',\gamma',\ip')$ on a nilpotent Lie group $N'$ of dimension $n$ is
isometric-isomorphic to a $\mu\in\nca_{\gamma}$.
\end{proposition}

\begin{proof}
We can assume that the Lie algebra of $N'$ is $(\ngo,\lambda)$ for some $\lambda\in\nca$.  There exist
$\vp\in\Gl(n)$ and $\psi\in\Or(\ngo,\ip)$ such that $\vp.\ip'=\ip$ and $\psi(\vp.\gamma')=\gamma$.  Thus the Lie
algebra isomorphism $\psi\vp:(\ngo,\lambda)\mapsto(\ngo,\mu)$, where $\mu=\psi\vp.\lambda$, satisfies
$\psi\vp.\ip'=\ip$ and $\psi\vp.\gamma'=\gamma$ and so it defines an isometry
$$
(N',\gamma',\ip')\mapsto(N_{\mu},\gamma,\ip)
$$
which is also an isomorphism between the class-$\gamma$ nilpotent Lie groups $(N',\gamma')$ and $(N,\gamma)$,
concluding the proof.
\end{proof}

According to the above proposition and identification (\ref{ideg}), the orbit $G_{\gamma}.\mu$ parameterizes all
the left invariant metrics which are compatible with $(N_{\mu},\gamma)$ and hence we may view $\nca_{\gamma}$ as
the space of all class-$\gamma$ metric structures on nilpotent Lie groups of dimension $n$. Since two metrics
$\mu,\lambda\in\nca_{\gamma}$ are isometric if and only if they live in the same $K_{\gamma}$-orbit, where
$K_{\gamma}=G_{\gamma}\cap\Or(\ngo,\ip)$ (see Appendix \ref{app}), we have that $\nca_{\gamma}/K_{\gamma}$
parameterizes class-$\gamma$ metric nilpotent Lie groups of dimension $n$ up to isometry and
$G_{\gamma}.\mu/K_{\gamma}$ do the same for all the compatible metrics on $(N_{\mu},\gamma)$.

We now recall a crucial fact which is the link between the study of left invariant compatible metrics for
geometric structures on nilpotent Lie groups and the results from real geometric invariant theory exposed in
Section \ref{git}.  The interplay is based on the identification given in (\ref{ideg}), and it will be used in
the proofs of the remaining results of this section.  The proof of the following proposition follows just from a
simple comparison between formulas (\ref{ricci}) and (\ref{defm2}).

\begin{proposition}\label{momric}
Let $m:V\mapsto\pg$ and $m_{\gamma}:V\mapsto\pg_{\gamma}$  be the moment maps for the actions of $\Gl(n)$ and
$G_{\gamma}$ on $V=\lam$, respectively {\rm (see Examples \ref{gln}, \ref{momg})}, where $\pg$ is the space of
symmetric maps of $(\ngo,\ip)$ and $\pg_{\gamma}$ the subspace of those leaving $\gamma$ invariant.
\begin{itemize}
\item[(i)] For each $\mu\in\nca\subset V$,
$$
m(\mu)=8\Ric_{\mu},
$$
where $\Ric_{\mu}$ is the Ricci operator of the Riemannian manifold $(N_{\mu},\ip)$.

\item[(ii)] For each $\mu\in\nca_{\gamma}\subset V$,
$$
m_{\gamma}(\mu)=8\Ricg_{\mu},
$$
where $\Ricg_{\mu}$ is the invariant Ricci operator of $(N_{\mu},\gamma,\ip)$, that is, the orthogonal
projection of the Ricci operator $\Ric_{\mu}$ on $\pg_{\gamma}$.
\end{itemize}
\end{proposition}

Part (i) will be applied in Section \ref{einstein} to consider the case $\gamma=0$, which is strongly related
with the study of Einstein solvmanifolds.

Let us now go back to our search for the best compatible metric.  The identification (\ref{ideg}) allows us to
view each point of the variety $\nca_{\gamma}$ as a class-$\gamma$ metric structure on a nilpotent Lie group of
dimension $n$. In this light, it is natural to consider the functional $F:\nca_{\gamma}\mapsto\RR$ given by
$F(\mu)=\tr(\Ricg_{\mu})^2$, which in some sense measures how far is the metric $\mu$ from having
$\Ricg_{\mu}=0$, which is the goal proposed in Section \ref{geometric} (see (\ref{uncondition})).  The critical
points of $F$ may be therefore considered compatible metrics of particular significance.  However, we should
consider some normalization since $F(t\mu)=t^4F(\mu)$ for all $t\in\RR$.

For any $\mu\in\nca_{\gamma}$ we have that $\scalar(\mu)=-\unc ||\mu||^2$ (see (\ref{ricci})), which says that
normalizing by scalar curvature and by the spheres of $V$ is equivalent:
$$
\{\mu\in\nca_{\gamma}:\scalar(\mu)=s\}=\{\mu\in\nca_{\gamma}:||\mu||^2=-4s\}, \quad \forall s<0.
$$
The critical points of $F:\PP V\mapsto\RR$, $F([\mu])=\tr(\Ricg_{\mu})^2/||\mu||^4$, which lie in
$\PP\nca_{\gamma}=\pi(\nca_{\gamma})$ appears then as very natural candidates, since it is like we are
restricting $F$ to the subset of all class-$\gamma$ metric structures having a given scalar curvature.

It follows from Proposition \ref{momric}, (ii), that
$$
F([\mu])=\frac{1}{64}||m_{\gamma}([\mu])||^2,
$$
where $m_{\gamma}:\PP V\mapsto\pg_{\gamma}$ is the moment map for the action of $G_{\gamma}$ on $\PP V$.  We
then obtain from Lemma \ref{marian2} and (\ref{actiong}) that
$$
\grad(F)_{[\mu]}=-\frac{1}{16}\pi^*\delta_{\mu}(\Ricg_{\mu}), \qquad ||\mu||=1,
$$
where $\pi^*:V\mapsto\tang_{[\mu]}\PP V$ denotes the derivative of the projection map $\pi:V\mapsto\PP V$.  This
shows that $[\mu]\in\PP V$ is a critical point of $F$ if and only if $\Ricg_{\mu}=cI+D$ for some $c\in\RR$ and
$D\in\Der(\mu)$ ($=\Ker{\delta_{\mu}}$).  By applying Theorem \ref{marian} to our situation we obtain the main
result of this section.

\begin{theorem}\label{equiv1g}
Let $F:\PP V\mapsto\RR$ be defined by $F([\mu])=\tr(\Ricg_{\mu})^2/||\mu||^4$.  Then for $\mu\in V$ the
following conditions are equivalent:
\begin{itemize}
\item[(i)] $[\mu]$ is a critical point of $F$.

\item[(ii)] $F|_{G_{\gamma}.[\mu]}$ attains its minimum value at $[\mu]$.

\item[(iii)] $\Ricg_{\mu}=cI+D$ for some $c\in\RR$, $D\in\Der(\mu)$.
\end{itemize}
Moreover, all the other critical points of $F$ in the orbit $G_{\gamma}.[\mu]$ lie in $K_{\gamma}.[\mu]$.
\end{theorem}

We now rewrite the above result in geometric terms, by using the identification (\ref{ideg}), Proposition
\ref{soliton1g} and Definition \ref{minimalg}.

\begin{theorem}\label{equiv2g}
Let $(N,\gamma)$ be a nilpotent Lie group endowed with an invariant geometric structure $\gamma$.  Then the
following conditions on a left invariant Riemannian metric $\ip$ which is compatible with $(N,\gamma)$ are
equivalent:
\begin{itemize}
\item[(i)] $\ip$ is minimal.

\item[(ii)] $\ip$ is an invariant Ricci soliton.

\item[(iii)] $\Ricg_{\ip}=cI+D$ for some $c\in\RR$, $D\in\Der(\ngo)$.
\end{itemize}
Moreover, there is at most one compatible left invariant metric on $(N,\gamma)$ up to isometry  (and scaling)
satisfying any of the above conditions.
\end{theorem}

Recall that the proof of this theorem does not use the integrability of $\gamma$, and so it is valid for the
`almost' versions as well.

We also note that part (i) of Theorem \ref{equiv1g} makes possibly the study of minimal compatible metrics by a
variational method.  Indeed, the projective algebraic variety $\PP\nca_{\gamma}$ may be viewed as the space of
all class-$\gamma$ metric structures on $n$-dimensional nilpotent Lie groups with a given scalar curvature, and
those which are minimal are precisely the critical points of $F:\PP\nca_{\gamma}\mapsto\RR$.  This variational
approach will be used quite often in the search for explicit examples in the following sections.

The above theorems propose then as privileged these compatible metrics called minimal, which have a neat
characterization (see (iii)), are critical points of a natural curvature functional (square norm of Ricci),
minimize such a functional when restricted to the compatible metrics for a given geometric structure, and are
solitons for a natural evolution flow.  Moreover, the uniqueness up to isometry of such special metrics holds.
But a remarkable weakness of this approach is the existence problem; the theorems do not even suggest when there
do exist such a distinguished metric.  How special are the symplectic or (almost-) complex structures admitting
a minimal metric?.  So far, we know how to deal with this `existence question' only by giving several examples,
which is the goal of Sections \ref{symp}-\ref{einstein}.  The only obstruction we have found is in the case
$\RR I\not\subset\ggo_{\gamma}$, namely when $\Der(\ngo)$ is nilpotent.  These Lie algebras are called {\it
characteristically nilpotent} and have been extensively studied in the last years, but we could not find any example
of a characteristically nilpotent Lie algebra admitting a symplectic structure.  Thus we have no any non-existence
example yet.

\begin{corollary}\label{isoiso}
Let $\gamma,\gamma'$ be two geometric structures on a nilpotent Lie group $N$, and assume that
they admit minimal compatible metrics $\ip$ and $\ip'$, respectively.  Then $\gamma$ is isomorphic to $\gamma'$
if and only if there exists $\vp\in\Aut(\ngo)$ and $c>0$ such that $\gamma'=\vp.\gamma$ and
$$
\la\vp X,\vp Y\ra'=c\la X,Y\ra \qquad \forall X,Y\in\ngo.
$$
In particular, if $\gamma$ and $\gamma'$ are isomorphic then their respective minimal compatible metrics are
necessarily isometric up to scaling (recall that $c=1$ when $\scalar(\ip)=\scalar(\ip')$).
\end{corollary}

We have here a very useful tool to distinguish two geometric structures.  Indeed, the corollary allows us to do
it by looking at their respective minimal compatible metrics, that is, with Riemannian data.  This is a
remarkable advantage since we suddenly have a great deal of invariants.  This method will be used in the subsequent
sections to find explicit continuous families depending on $1$, $2$ and $3$ parameters of pairwise
non-isomorphic geometric structures in low dimensions, mainly by using only one Riemannian invariant: the
eigenvalues of the Ricci operator.

\begin{remark}\label{liegroups}
{\rm
The Ricci curvature operator of a left invariant metric $\ip$ on a Lie group $G$ is given by
$$
\Ric_{\ip}=R_{\ip}-\unm B_{\ip}-D_{\ip},
$$
where $R_{\ip}$ is defined by (\ref{ricci}), $B_{\ip}$ is the Killing form of the Lie algebra $\ggo$ of $G$ in
terms of $\ip$, $D_{\ip}$ is the symmetric part of $\ad{Z_{\ip}}$ and $Z_{\ip}\in\ggo$ is defined by $\la
Z_{\ip},X\ra=\tr(\ad{X})$ for any $X\in\ggo$.  Recall that $Z_{\ip}=0$ if and only if $\ggo$ is unimodular, and
that $\Ric_{\ip}=R_{\ip}$ in the nilpotent case.  If we consider the tensor $R$ instead of the Ricci tensor, and
the variety of all Lie algebras $\lca$ rather than just the nilpotent ones, to define and state all the notions,
flows, identifications and results in Sections \ref{geometric} and \ref{var}, then everything is still valid for
Lie groups in general, with the only exception of the first part of Corollary \ref{isoiso}. We only have to
consider the corresponding invariant part $R^{\gamma}$ and replace $\scalar(\ip)$ with $\tr{R_{\ip}}$ each time it
appears.  The only detail to be careful with is that if two $\mu,\lambda\in\lca_{\gamma}$ lie in the same
$K_{\gamma}$-orbit then they are isometric, but the converse might not be true.  Recall that the uniqueness
result in Theorem \ref{equiv2g} is nevertheless valid.

The reason why we decided to work only in the nilpotent case is that, at least at first sight, the use of this
`unnamed' tensor $R$ make minimal and soliton metrics, as well as the functionals and evolution flows, into
concepts lacking in geometric sense.  For instance, we have found ourselves with the unpleasant fact that some
K$\ddot{{\rm a}}$hler metrics on solvable Lie groups would not be minimal viewed as compatible metrics for the
corresponding symplectic structures, in spite of $\Ricac_{\ip}=0$.

For a compact simple Lie group, $-B$ is minimal for the case $\gamma=0$, and if $\ggo=\kg\oplus\pg$ is the
Cartan decomposition for a non-compact semi-simple Lie algebra $\ggo$, then it is easy to see that the metric
$\ip$ given by $\la\kg,\pg\ra=0$, $\ip|_{\kg\times\kg}=-B$ and $\ip|_{\pg\times\pg}=B$ is minimal as well.    }
\end{remark}

\section{Symplectic structures}\label{symp}

\subsection{Metrics with hermitian Ricci tensor and the anti-complexified Ricci flow}\label{hermricci}  Let
$(M,\omega)$ be a symplectic
manifold, that is, a differentiable manifold $M$ endowed with a global $2$-form $\omega$ which is closed
($\dif\omega=0$) and non-degenerate ($\omega^n\ne 0$).  A Riemannian metric $g$ on $M$ is said to be {\it
compatible} with $\omega$ if there exists an almost-complex structure $J_g$ (i.e. a $(1,1)$-tensor field with
$J_g^2=-I$) such that
$$
\omega=g(\cdot,J_g\cdot).
$$
In that case $J_g$ is uniquely determined by $g$, and one may also define that an almost-complex structure $J$
is {\it compatible} with $\omega$ if
$$
g_J=\omega(\cdot,J\cdot)
$$
determines a Riemannian metric, which is again uniquely determined by $J$.  In such a way we are really talking
about compatible pairs $(g,J)=(g,J_g)=(g_J,J)$, and the triple $(\omega,g,J_g)$ is called an {\it
almost-K$\ddot{{\rm a}}$hler} structure on $M$.

It is well known that for any symplectic manifold there always exist a compatible metric.  Moreover, the space
$\cca=\cca(M,\omega)$ of all compatible metrics is usually huge; recall for instance that the group of all
symplectomorphisms (i.e. diffeomorphisms $\vp$ of $M$ such that $\vp^*\omega= \omega$) acts on $\cca$.

We fix from now on a symplectic manifold $(M,\omega)$.  Let $\Ric_g$ and $\nabla_g$ denote the Ricci operator
and the Levi-Civita conexion of a compatible metric $g$, respectively.  The most famous conditions to ask g to
satisfy are Einstein (i.e. $\Ric_g=cI$) and K$\ddot{{\rm a}}$hler (i.e. $\nabla_gJ_g =0$), which are both very
strong and share the following property.

\begin{definition}\label{hrt}
{\rm We say that $g$ has {\it hermitian} (or $J$-{\it invariant}) {\it Ricci tensor} or that $J_g$ is {\it
harmonic}, if $\Ric_gJ_g= J_g\Ric_g$.}
\end{definition}
Examples of compatible metrics with hermitian Ricci tensor which are neither Einstein nor K$\ddot{{\rm a}}$hler
are known in any dimension $2n\geq 6$ (see \cite{ApsGdc, DvdMsk,LeWng}).  We will show in Section \ref{nilpsymp}
that a symplectic nilpotent Lie group can never admit a compatible metric with hermitian Ricci tensor unless it
is abelian.

A classical approach to searching for distinguished metrics is the variational one, that is, to consider
critical points of natural functionals of the curvature on the space of all metrics of a given class.  For
instance, if $M$ is compact, D. Hilbert \cite{Hlb} proved that Einstein metrics on $M$ are precisely the
critical points of the total scalar curvature functional
$$
S:\mca_1\mapsto\RR, \qquad S(g)=\int_M\scalar(g)\dif\nu_g,
$$
where $\mca_1$ is the space of all Riemannian metrics on $M$ with volume equal to $1$.  Since the set of
compatible metrics $\cca$ is smaller, one should expect a weaker critical point condition for
$S:\cca\mapsto\RR$.  Another natural functional in our setup would be
$$
K:\cca\mapsto\RR, \qquad K(g)=\int_M||\nabla_gJ_g||^2\dif\nu_g,
$$
for which K$\ddot{{\rm a}}$hler metrics are precisely the global minima. D. Blair and S. Ianus proved that,
curiously enough, both functionals $S$ and $K$ have the same critical points on $\cca$.

\begin{theorem}\label{critherm} \cite{BlrIns}
Let $(M,\omega)$ be a symplectic manifold and $\cca$ the set of all compatible metrics.  Then $g\in\cca$ is a
critical point of $S:\cca\mapsto\RR$ or $K:\cca\mapsto\RR$ if and only if $g$ has hermitian Ricci tensor.
\end{theorem}

This result and the above considerations do suggest that the compatible metrics with hermitian Ricci tensor (if
any) are really `good friends' of the symplectic structure.

In \cite{LeWng}, H-V Le and G. Wang approach the problem of the existence of such metrics by considering an
evolution flow inspired in the Ricci flow introduced by R. Hamilton \cite{Hml1}.  If $\ricci_g$ is the Ricci
tensor of a compatible metric $g$, then consider the orthogonal decomposition
\begin{equation}\label{ac-c}
\ricci_g=\ricciac_g+\riccic_g,
\end{equation}
where $\ricciac_g=\unm(\ricci_g-\ricci_g(J_g\cdot,J_g\cdot))$ and
$\riccic_g=\unm(\ricci_g+\ricci_g(J_g\cdot,J_g\cdot))$ are the {\it anti-complexified} and {\it complexified}
parts of $\ricci_g$, respectively. In this way, $g$ has hermitian Ricci tensor if and only if $\ricciac_g=0$,
and since the gradient of the functional $K$ equals $-\ricciac_g$ it is natural to consider the negative
gradient flow equation
\begin{equation}\label{acrf}
\ddt g(t)=\ricciac_{g(t)},
\end{equation}
for a curve $g(t)$ of metrics, which is called in \cite{LeWng} the {\it anti-complexified Ricci flow}.  Recall
that the fixed points of (\ref{acrf}) are precisely the metrics with hermitian Ricci tensor.  The main result in
\cite{LeWng} is the short time existence and uniqueness of the solution to (\ref{acrf}) when $M$ is compact.

\subsection{Symplectic nilpotent Lie groups}\label{nilpsymp}
Let $N$ be a real $2n$-dimensional nilpotent Lie group with Lie algebra $\ngo$, whose Lie bracket is denoted by
$\mu :\ngo\times\ngo\mapsto\ngo$.  An invariant {\it symplectic} structure on $N$ is defined by a $2$-form
$\omega$ on $\ngo$ satisfying
$$
\omega(X,\cdot)\equiv 0 \quad \mbox{if and only if} \quad X=0 \quad (\mbox{non-degenerate}),
$$
and for all $X,Y,Z\in\ngo$
\begin{equation}\label{closed}
\omega(\mu(X,Y),Z)+\omega(\mu(Y,Z),X)+\omega(\mu(Z,X),Y)=0 \quad ({\rm closed}, \dif\omega=0).
\end{equation}
Fix a symplectic nilpotent Lie group $(N,\omega)$.  A left invariant Riemannian metric which is compatible with
$(N,\omega)$ is determined by an inner product $\ip$ on $\ngo$ such that if
$$
\omega(X,Y)=\la X,J_{\ip}Y\ra\quad\forall\; X,Y\in\ngo\quad{\rm then}\quad J_{\ip}^2=-I.
$$
For the geometric structure $\gamma=\omega$ we have that
$$
G_{\gamma}=\Spe(n,\RR)=\{ g\in\Gl(2n):g^tJg=J\}, \qquad K_{\gamma}=\U(n),
$$
and the Cartan decomposition of $\ggo_{\gamma}=\spg(n,\RR)=\{ A\in\glg(2n):A^tJ+JA=0\}$ is given by
$$
\spg(n,\RR)=\ug(n)\oplus\pg_{\gamma}, \qquad \pg_{\gamma}=\{ A\in\pg:AJ=-JA\}.
$$
Thus the invariant Ricci tensor $\riccig$ coincides with the anti-complexified Ricci tensor (see (\ref{ac-c}))
and for any $\ip\in\cca$,
\begin{equation}\label{ac-cop}
\Ricg_{\ip}=\Ricac_{\ip}=\unm(\Ric_{\ip}+J_{\ip}\Ric_{\ip}J_{\ip}).
\end{equation}
This implies that our `goal' condition $\Ricg_{\ip}=0$ (see (\ref{uncondition})) is equivalent to have hermitian
Ricci tensor. Also, the evolution flow considered in Section \ref{geometric} is not other than the
anti-complexified Ricci flow.

Concerning the search for the best compatible left invariant metric for a symplectic nilpotent Lie group, and in
view of the facts exposed in Section \ref{hermricci}, our first result is negative.

\begin{proposition}\label{noherm}
Let $(N,\omega)$ be a symplectic nilpotent Lie group.  Then $(N,\omega)$ does not admit any compatible left
invariant metric with hermitian Ricci tensor, unless $N$ is abelian.
\end{proposition}

\begin{proof}
We first note that since $\mu$ is nilpotent the center $\zg$ of $(\ngo,\mu)$ meets non-trivially the derived Lie
subalgebra $\mu(\ngo,\ngo)$, unless $\mu=0$ (i.e. $\ngo$ abelian).  Assume that $\ip\in\cca(N,\omega)$ has
hermitian Ricci tensor and consider the orthogonal decomposition $\ngo=\vg\oplus\mu(\ngo,\ngo)$.  If $Z\in\zg$
then $JZ\in\vg$. In fact, it follows from (\ref{closed}) that
$$
\la\mu(X,Y),JZ\ra=\omega(\mu(X,Y),Z)=0 \qquad \forall\; X,Y\in\ngo.
$$
Now, the above equation, the fact that $\Ric_{\ip}J=J\Ric_{\ip}$ and the definition of $\Ric_{\ip}$ (see
(\ref{ricci})) imply that
$$
0\leq\la\Ric_{\ip}Z,Z\ra=\la\Ric_{\ip}JZ,JZ\ra\leq 0,
$$
and hence
$$
\unc\sum_{ij}\la\mu(X_i,X_j),Z\ra^2=\la\Ric_{\ip}Z,Z\ra=0 \qquad \forall\; Z\in\zg.
$$
Thus $\mu(\ngo,\ngo)\perp\zg$ and so $\ngo$ must be abelian by the observation made in the beginning of the
proof.
\end{proof}

We now review the variational approach developed in Section \ref{var} and obtain some applications.  Fix a
non-degenerate $2$-form $\omega$ on $\ngo$, and let $\Spe(n,\RR)$ denote the subgroup of $\Gl(2n)$ preserving
$\omega$, that is,
$$
\Spe(n,\RR)=\{ \vp\in\Gl(2n):\omega(\vp X,\vp Y)=\omega(X,Y) \quad\forall\; X,Y\in\ngo\}.
$$
Consider the algebraic subvariety $\nca_s:=\nca_{\gamma}\subset\nca$ given by
$$
\nca_s=\{\mu\in\nca:\dif_{\mu}\omega=0\},
$$
that is, those nilpotent Lie brackets for which $\omega$ is closed (see (\ref{closed})).  By fixing an inner
product $\ip$ on $\ngo$ satisfying that
$$
\omega=\la\cdot,J\cdot\ra\qquad{\rm with}\qquad J^2=-I,
$$
(\ref{ideg}) identify each $\mu\in\nca_s$ with the almost-K$\ddot{{\rm a}}$hler manifold
$(N_{\mu},\omega,\ip,J)$.  The action of $\Spe(n,\RR)$ on $\nca_s$ has the following interpretation:  each
$\vp\in\Spe(n,\RR)$ determines a Riemannian isometry which is also a symplectomorphism
$$
(N_{\vp.\mu},\omega,\ip,J)\mapsto (N_{\mu},\omega,\la\vp\cdot,\vp\cdot\ra,\vp^{-1}J\vp)
$$
by exponentiating the Lie algebra isomorphism $\vp^{-1}:(\ngo,\vp.\mu)\mapsto(\ngo,\mu)$.

We have seen in Proposition \ref{noherm} that the hermitian Ricci condition $\Ricac=0$ is a very nice but
forbidden condition for a symplectic non-abelian nilpotent Lie group. However, we can get as close as we want.

\begin{proposition}\label{epsilon}
Let $(N,\omega)$ be a symplectic nilpotent Lie group.  Then for any $\epsilon>0$ there exists a compatible left
invariant metric $\ip$ such that $\tr(\Ricac_{\ip})^2<\epsilon$.
\end{proposition}

\begin{proof}
In view of the identification (\ref{ideg}), this is equivalent to prove that for any $\mu\in\nca_s$ there exists
$\lambda\in\Spe(n,\RR).\mu$ such that $\tr(\Ricac_{\lambda})^2<\epsilon$.  We have that an orbit
$\Spe(n,\RR).\mu$ is closed if and only if $\Ricac_{\lambda}=0$ for some $\lambda\in\Spe(n,\RR).\mu$ (see
Theorem \ref{RS}, (iii) and Proposition \ref{momric}, (ii)).  It then follows from Proposition \ref{noherm} that
for $\mu\in\nca_s$ the orbit $\Spe(n,\RR).\mu$ can never be closed, unless $\mu=0$.  Thus
$0\in\overline{\Spe(n,\RR).\mu}$ for any $\mu\in\nca_s$ (see Theorem \ref{RS}, (v)), and the continuity of the
function $\lambda\mapsto\tr(\Ricac_{\lambda})^2$ concludes the proof.
\end{proof}

\begin{proposition}\label{sc1}
Let $(N,\omega)$ be a symplectic nilpotent Lie group.  For any real number $s<0$ there exists a compatible
metric $\ip$ such that $\scalar(\ip)=s$.
\end{proposition}

\begin{proof}
By considering the identification (\ref{ideg}) and taking in account that $\scalar(\mu)=-\unc||\mu||^2$ (see
(\ref{ricci})), we get that this proposition is equivalent to the following fact: for each $\mu\in\nca_s$ the
orbit $\Spe(n,\RR).\mu$ meets all the spheres of $V$.  In the proof of Proposition \ref{epsilon} we have seen
that $0\in\overline{\Spe(n,\RR).\mu}$ for any $\mu\in\nca_s$, and so by Theorem \ref{RS}, (iv), there exists
$A\in\pg$ such that $\lim_{t\to-\infty}\exp{tA}.\mu=0$.

If $f(t)=||\exp{tA}.\mu||^2$ then $f''(t)>0$ for all $t\in\RR$ (see \cite[Lemma 3.1]{RchSld}) and so
$\lim_{t\to\infty}f(t)=+\infty$, concluding the proof.
\end{proof}

\subsection{Examples}
Let $\ngo$ be a $2n$-dimensional vector space with basis $\{ X_1,...,X_{2n}\}$ over $\RR$, and consider the
non-degenerate $2$-form
$$
\omega=\alpha_1\wedge\alpha_{2n}+...+\alpha_n\wedge\alpha_{n+1},
$$
where $\{\alpha_1,...,\alpha_{2n}\}$ is the dual basis of $\{ X_i\}$.  For the compatible inner product $\la
X_i,X_j\ra=\delta_{ij}$ we have that $\omega=\la\cdot,J\cdot\ra$ for
$$
J=\left[\begin{smallmatrix} &&&&&-1\\ &0&&&\cdot&\\ &&&-1&&\\ &&1&&&\\ &\cdot&&&0&\\
1&&&&&\end{smallmatrix}\right].
$$
In all the examples the symplectic structure will be $\omega$, the almost-complex structure $J$ and the
compatible metric $\ip$.  We will vary Lie brackets and use constantly identification (\ref{ideg}).

\begin{example}\label{heis}
{\rm Let $\mu_n$ the $2n$-dimensional Lie algebra whose only non-zero bracket is
$$
\mu_n(X_1,X_2)=X_3,
$$
that is, $\mu_n$ is isomorphic to $\hg_3\oplus\RR^n$, where $\hg_3$ is the $3$-dimensional Heisenberg Lie
algebra.  Recall that $(N_{\mu_2},\omega)$ is precisely the simply connected cover of the famous
Kodaira-Thurston manifold.  It is easy to prove that $\overline{\Spe(n,\RR).\mu_n}=\Spe(n,\RR).\mu_n\cup\{ 0\}$,
and so the orbit $\Spe(n,\RR).[\mu_n]$ is closed in $\PP V$.  This implies that the functional $F$ from Theorem
\ref{equiv1g} must attain its minimum value on $\Spe(n,\RR).[\mu_n]$ and hence there exists a metric compatible
with $(N_{\mu_n},\omega)$ which is minimal.  In fact, the inner product $\la X_i,X_j\ra=\delta_{ij}$ satisfies
$$
\begin{array}{l}
\Ricac_{\ip}=-\unc \left[\begin{smallmatrix} 1&&&&&&&&\\ &1&&&&&&&\\ &&-1&&&&&&\\ &&&0&&&&&\\ &&&&\cdot&&&&\\
&&&&&0&&&\\ &&&&&&1&&\\ &&&&&&&-1&\\ &&&&&&&&-1\end{smallmatrix}\right] = -\frac{3}{4}I+\unc
\left[\begin{smallmatrix} 2&&&&&&&&\\ &2&&&&&&&\\ &&4&&&&&&\\ &&&3&&&&&\\
&&&&\cdot&&&&\\ &&&&&3&&&\\ &&&&&&2&&\\ &&&&&&&4&\\ &&&&&&&&4\end{smallmatrix}\right], \qquad 2n\geq 8, \\

\Ricac_{\ip}=-\unc \left[\begin{smallmatrix} 1&&&&&\\ &1&&&&\\ &&-1&&&\\ &&&1&&\\ &&&&-1&\\
&&&&&-1\end{smallmatrix}\right] = -\frac{3}{4}I+\unc \left[\begin{smallmatrix} 2&&&&&\\ &2&&&&\\ &&4&&&\\ &&&2&&\\
&&&&4&\\ &&&&&4\end{smallmatrix}\right], \qquad 2n=6, \\

\Ricac_{\ip}=-\unc \left[\begin{smallmatrix} 1&&&\\ &2&&\\ &&-2&\\ &&&-1\end{smallmatrix}\right] =
-\frac{5}{4}I+\unc \left[\begin{smallmatrix} 4&&&\\ &3&&\\ &&7&\\ &&&6\end{smallmatrix}\right], \qquad 2n=4,
\end{array}
$$
and hence $\Ricac_{\ip}\in\RR I+\Der(\mu_n)$ in all the cases.  Moreover, it follows from the closeness of
$\Spe(n,\RR).[\mu_n]$ that $F$ must also attain its maximum value, and therefore
$\Spe(n,\RR).[\mu_n]=\U(n).[\mu_n]$ by uniqueness in Theorem \ref{equiv1g}.  This implies that there is only one
left invariant metric compatible with $(N_{\mu_n},\omega)$ up to isometry, often called the Abbena metric in the
case $n=2$. }
\end{example}

\begin{example}\label{otra4}
{\rm Consider the $4$-dimensional Lie algebra given by
$$
\lambda(X_1,X_2)=X_3, \qquad \lambda(X_1,X_3)=X_4.
$$
The compatible metric $\la X_i,X_j\ra=\delta_{ij}$ is minimal for $(N_{\lambda},\omega)$ since
$$
\Ricac_{\ip}=-\unc \left[\begin{smallmatrix} 3&&&\\ &1&&\\ &&-1&\\ &&&-3\end{smallmatrix}\right] =
-\frac{5}{4}I+\unm\left[\begin{smallmatrix} 1&&&\\ &2&&\\ &&3&\\ &&&4\end{smallmatrix}\right]\in\RR
I+\Der(\lambda).
$$
  }
\end{example}

It is well-known that $(N_{\mu_2},\omega)$ and $(N_{\lambda},\omega)$ are the only symplectic nilpotent Lie
groups in dimension $4$, and then the existence of minimal compatible metrics in the case $2n=4$ follows.

\begin{example}\label{abc}
{\rm Let $\mu=\mu(a,b,c)$ be the $6$-dimensional $2$-step nilpotent Lie algebra defined by
$$
\mu(X_1,X_2)=aX_4, \quad \mu(X_1,X_3)=bX_5, \quad \mu(X_2,X_3)=cX_6.
$$
It is easy to check that $\mu\in\nca_s$ if and only if $a-b+c=0$.  We can also get from a simple calculation
that
$$
\begin{array}{rl}
\Ricac_{\mu}&=-\unc (a^2+b^2+c^2) \left[\begin{smallmatrix} 1&&&&&\\ &1&&&&\\ &&1&&&\\
&&&-1&&\\ &&&&-1&\\ &&&&&-1\end{smallmatrix}\right] \\
&=-\frac{3}{4}(a^2+b^2+c^2)I+\unm (a^2+b^2+c^2)\left[\begin{smallmatrix} 1&&&&&\\ &1&&&&\\ &&1&&&\\
&&&2&&\\ &&&&2&\\ &&&&&2 \end{smallmatrix}\right]\in\RR I+\Der(\mu),
\end{array}
$$
and so the whole family $\{[\mu(a,b,c)]:a-b+c=0\}\subset\PP\nca_s$ consists of critical points of $F$.  We
assume that $a^2+b^2+c^2=2$ in order to avoid homothetical changes, which is equivalent to $\scalar(\mu)=-1$.
The Ricci operator on the center $\zg=\la X_4,X_5,X_6\ra_{\RR}$ is given by
$$
\Ric_{\mu}=\unm\left[\begin{smallmatrix} a^2&&\\ &b^2&\\ &&c^2\end{smallmatrix}\right],
$$
and thus the curve $\{\mu_{st}=\mu(s,s+t,t):s^2+st+t^2=1,\; 0\leq t\leq\frac{1}{\sqrt{3}}\}$ is pairwise
non-isometric.  It then follows from Corollary \ref{isoiso} that $(N_{\mu_{st}},\omega)$ is a curve of pairwise
non-isomorphic symplectic nilpotent Lie groups.  In terms of the notation in \cite[Table A.1]{Slm}, we have that
$\mu_{01}\simeq (0,0,0,0,12,13)$ and $\mu_{st}\simeq (0,0,0,12,13,23)$ for any $0<t\leq\frac{1}{\sqrt{3}}$.  We
note that this curve coincides with the curve of pairwise non-isomorphic symplectic structures denoted by
$\omega_1(t)$ in \cite[Theorem 3.1, 18]{KhkGzMdn}, and then this example shows that any symplectic structure in
such a curve admits a compatible metric which is minimal  }
\end{example}

\begin{example}\label{m26}
{\rm We now take advantage of the variational nature of Theorem \ref{equiv1g} to find explicit examples of
minimal compatible metrics.  Consider for each $6$-upla $\{ a,...,f\}$ of real numbers the skew-symmetric
bilinear form $\mu=\mu(a,b,c,d,e,f)\in V=\Lambda^2(\RR^6)^*\otimes\RR^6$ defined by
$$
\begin{array}{l}
\mu(X_1,X_2)=aX_3, \qquad\mu(X_1,X_3)=bX_4, \qquad \mu(X_1,X_4)=cX_5,\\

\mu(X_1,X_5)=dX_6, \qquad\mu(X_2,X_3)=eX_5,  \qquad\mu(X_2,X_4)=fX_6.
\end{array}
$$

Our plan is to find first the critical points of $F$ restricted to the set $\{[\mu(a,...,f)]:a,...,f\in\RR\}$
and after that to show by using the characterization given in part (iii) of the theorem that they are really
critical points of $F:\PP V\mapsto\RR$.  We can see by a simple computation that $\Ricac_{\mu}$ is given by the
diagonal matrix with entries
$$
\begin{array}{rl}
\Ricac_{\mu}&=-\unc\left[\begin{smallmatrix} a^2+b^2+c^2+2d^2+f^2\\ a^2+c^2-d^2+2e^2+f^2\\
-a^2+2b^2-c^2+e^2-f^2\\
a^2-2b^2+c^2-e^2+f^2\\ -a^2-c^2+d^2-2e^2-f^2\\ -a^2-b^2-c^2-2d^2-f^2\end{smallmatrix}\right],
\end{array}
$$
and hence we are interested in the critical points of
$$
\begin{array}{rl}
F(\mu)=&\tr(\Ricac_{\mu})^2=F(a,...,f)\\
=&\frac{1}{8}\Big((a^2+b^2+c^2+2d^2+f^2)^2+(a^2+c^2-d^2+2e^2+f^2)^2 \\
&+(-a^2+2b^2-c^2+e^2-f^2)^2\Big)
\end{array}
$$
restricted to any leaf of the form $a^2+b^2+c^2+d^2+e^2+f^2\equiv$ const., which are easily seen to depend of
three parameters.  We still have to impose the Jacobi and closeness conditions on these critical points (or
equivalently to find the intersection with $\nca_s$), after which we obtain the following ellipse of symplectic
structures:
$$
\{\mu_{xy}=\mu(x,1,x+y,1,1,y):x^2+y^2+xy=1\}.
$$
It follows from the formula for $\Ricac_{\mu}$ given above that
$$
\Ricac_{\mu_{xy}}=-\unc\left[\begin{smallmatrix} 5&&&&&\\ &3&&&&\\ &&1&&&\\
&&&-1&&\\ &&&&-3&\\ &&&&&-5\end{smallmatrix}\right]
=-\frac{7}{4}I+\unm\left[\begin{smallmatrix} 1&&&&&\\ &2&&&&\\ &&3&&&\\
&&&4&&\\ &&&&5&\\ &&&&&6 \end{smallmatrix}\right]\in\RR I+\Der(\mu_{xy}),
$$
showing definitely that this is a curve of minimal compatible metrics.  We furthermore have that the Ricci
tensor of the metrics $\mu_{xy}$ is given by
$$
\Ric_{\mu_{xy}}=-\unc\left[\begin{smallmatrix} 4-y^2&&&&&\\ &2-xy&&&&\\ &&2-x^2&&&\\ &&&1-x^2&&\\ &&&&-1-xy&\\
&&&&&-1-y^2\end{smallmatrix}\right],
$$
which clearly shows that they are pairwise non-isometric for $x,y\geq 0$.  It then follows from Corollary
\ref{isoiso} that
$$
\{(N_{\mu_{xy}},\omega):x^2+y^2+xy=1,\; x,y\geq 0\}
$$
is a curve of pairwise non-isomorphic symplectic nilpotent Lie groups.  There are three $6$-dimensional
nilpotent Lie groups involved, $N_{\mu_{10}}$, $N_{\mu_{01}}$ and $N_{\mu_{xy}}$, $x,y>0$, denoted in
\cite[Table A.1]{Slm} by $(0,0,12,13,14,23+15)$, $(0,0,0,12,14-23,15+34)$ and $(0,0,12,13,14+23,24+15)$,
respectively.  }
\end{example}

\section{Complex structures}\label{complex}

Let $N$ be a real $2n$-dimensional nilpotent Lie group with Lie algebra $\ngo$, whose Lie bracket is denoted by
$\mu :\ngo\times\ngo\mapsto\ngo$.  An invariant {\it almost-complex} structure on $N$ is defined by a map
$J:\ngo\mapsto\ngo$ satisfying $J^2=-I$.  If in addition $J$ satisfies the integrability condition
\begin{equation}\label{integral}
\mu(JX,JY)=\mu(X,Y)+J\mu(JX,Y)+J\mu(X,JY), \qquad \forall X,Y\in\ngo,
\end{equation}
then $J$ is said to be a {\it complex} structure.

Fix an almost-complex nilpotent Lie group $(N,J)$.  A left invariant Riemannian metric which is {\it compatible}
with $(N,J)$, also called an {\it almost-hermitian metric}, is given by an inner product $\ip$ on $\ngo$ such
that
$$
\la JX,JY\ra=\la X,Y\ra \qquad \forall X,Y\in\ngo.
$$
We have for this particular geometric structure $\gamma=J$ that
$$
G_{\gamma}=\Gl(n,\CC)=\{ \vp\in\Gl(2n):\vp J=J\vp\}, \qquad K_{\gamma}=\U(n),
$$
and the Cartan decomposition of $\ggo_{\gamma}=\glg(n,\CC)=\{ A\in\glg(2n):AJ=JA\}$ is given by
$$
\glg(n,\CC)=\ug(n)\oplus\pg_{\gamma}, \qquad \pg_{\gamma}=\{ A\in\pg:AJ=JA\}.
$$
The invariant Ricci operator is then given by the complexified Ricci operator
$$
\Ricg_{\ip}=\Ricc_{\ip}=\unm(\Ric_{\ip}-J\Ric_{\ip}J)
$$
(see (\ref{ac-c})).  In this way, condition $\Ricg_{\ip}=0$ is equivalent to the Ricci operator anti-commute
with $J$. We do not know if this property has any special significance in complex geometry, but for instance it
holds for a K$\ddot{{\rm a}}$hler metric if and only if the metric is Ricci flat.  Anyway, as in the symplectic
case, the condition $\Ricc_{\ip}=0$ is also forbidden for non-abelian $N$ since
$\tr{\Ricc_{\ip}}=\scalar(\ip)<0$.

We now fix a map $J:\ngo\mapsto\ngo$ satisfying $J^2=-I$ and consider the algebraic subvariety
$\nca_c:=\nca_{\gamma}\subset\nca$ given by
$$
\nca_c=\{\mu\in\nca: \mbox{(\ref{integral}) holds}\},
$$
that is, those nilpotent Lie brackets for which $J$ is integrable and so define a complex structure on
$N_{\mu}$, the simply connected nilpotent Lie group with Lie algebra $(\ngo,\mu)$.

Fix also an inner product $\ip$ on $\ngo$ compatible with $J$, then (\ref{ideg}) identifies each $\mu\in\nca_c$
(or $\nca$) with the hermitian (or almost-hermitian) manifold $(N_{\mu},J,\ip)$.  If we use the same triple
$(\omega,J,\ip)$ to define and identify $\nca_s$ (see Section \ref{nilpsymp}) and $\nca_c$, then the
intersection of these varieties is $\nca_s\cap\nca_c=\{ 0\}$ since no non-abelian nilpotent Lie group can admit
a K$\ddot{{\rm a}}$hler metric.

We now give some examples.

\begin{example}\label{heisc}
{\rm Let $\mu_n$ be the $2n$-dimensional Lie algebra considered in Example \ref{heis}.  It is easy to check that
$\ip$ is also minimal as a compatible metric for the almost-complex nilpotent Lie group $(N_{\mu_n},J)$.  For
the $4$-dimensional Lie algebra in Example \ref{otra4}, we have that $\Ricc_{\ip}=-\unc I$ and hence this metric
is minimal for the almost-complex nilpotent Lie group $(N_{\lambda},J)$ as well.  }
\end{example}

For $\ngo_1=\RR^4$ and $\ngo_2=\RR^2$, consider the vector space $W=\Lambda^2\ngo_1^*\otimes\ngo_2$ of all
skew-symmetric bilinear maps $\mu:\ngo_1\times\ngo_1\mapsto\ngo_2$.  Any $6$-dimensional $2$-step nilpotent Lie
algebra with $\dim{\mu(\ngo,\ngo)}\leq 2$ can be modelled in this way.  Fix basis $\{ X_1,...,X_4\}$ and $\{
Z_1,Z_2\}$ of $\ngo_1$ and $\ngo_2$, respectively.  Each element in  $W$ will be described as
$\mu=\mu(a_1,a_2,...,f_1,f_2)$, where {\small
$$
\begin{array}{lll}
\mu(X_1,X_2)=a_1Z_1+a_2Z_2, & \mu(X_1,X_4)=c_1Z_1+c_2Z_2, & \mu(X_2,X_4)=e_1Z_1+e_2Z_2,\\
\mu(X_1,X_3)=b_1Z_1+b_2Z_2, & \mu(X_2,X_3)=d_1Z_1+d_2Z_2, & \mu(X_3,X_4)=f_1Z_1+f_2Z_2.
\end{array}
$$}
The complex structure and the compatible metric will be always defined by
$$
\begin{array}{lll}
J=\left[\begin{smallmatrix} 0&-1&&&&\\ 1&0&&&&\\ &&0&-1&&\\ &&1&0&&\\ &&&&0&-1\\ &&&&1&0
\end{smallmatrix}\right], && \la X_i,X_j\ra=\la Z_i,Z_j\ra=\delta_{ij}.
\end{array}
$$
If $A=(a_1,a_2)$, ... , $F=(f_1,f_2)$ and $JA=(-a_2,a_1)$, ... , $JF=(-f_2,f_1)$, then it is easy to check that
$J$ is integrable on $N_{\mu}$ (or $(N_{\mu},J)$ is a complex nilpotent Lie group) if and only if
\begin{equation}\label{integrable}
E=B+JD+JC,
\end{equation}
$J$ is {\it bi-invariant} (i.e. $\mu(JX,Y)=J\mu(X,Y)$) if and only if
\begin{equation}\label{biinv}
A=F=0, \qquad C=D=JB, \qquad E=-B,
\end{equation}
and $J$ is {\it abelian} (i.e. $\mu(JX,JY)=\mu(X,Y)$) if and only if
\begin{equation}\label{abelian}
E=B, \qquad D=-C.
\end{equation}
We note that the above conditions determine $\Gl(2,\CC)\times\Gl(1,\CC)$-invariant linear subspaces of $W$ of
dimensions $10$, $2$ and $8$, respectively.  For each $\mu\in W$, it follows from (\ref{ricci}) that the Ricci
operator of the almost-hermitian manifold $(N_{\mu},J,\ip)$ (see identification (\ref{ideg})) restricted to
$\ngo_1$, $\Ric_{\mu}|_{\ngo_1}$, is given by
\begin{equation}\label{ricciw}
-\unm\left[\begin{smallmatrix}
||A||^2+||B||^2+||C||^2 & \la B,D\ra+\la C,E\ra & -\la A,D\ra+\la C,F\ra & -\la A,E\ra-\la B,F\ra\\
\la B,D\ra+\la C,E\ra & ||A||^2+||D||^2+||E||^2 & \la A,B\ra+\la E,F\ra & \la A,C\ra-\la D,F\ra\\
-\la A,D\ra+\la C,F\ra & \la A,B\ra+\la E,F\ra & ||B||^2+||D||^2+||F||^2 & \la B,C\ra+\la D,E\ra\\
-\la A,E\ra-\la B,F\ra & \la A,C\ra-\la D,F\ra & \la B,C\ra+\la D,E\ra & ||C||^2+||E||^2+||F||^2
\end{smallmatrix}\right]
\end{equation}
and
$$
\Ric_{\mu}|_{\ngo_2}=\unm\left[\begin{smallmatrix} ||v_1||^2 & \la v_1,v_2\ra\\ \la v_1,v_2\ra & ||v_2||^2
\end{smallmatrix}\right], \quad v_i=(a_i,b_i,c_i,d_i,e_i,f_i),\quad i=1,2.
$$
Recall that if the complexified Ricci operator satisfied $\Ricc_{\mu}|_{\ngo_1}=pI$, $p\in\RR$, then $\mu$ is
minimal.  Indeed, since always $\Ricc_{\mu}|_{\ngo_2}=qI$ for some $q\in\RR$, we would have that
\begin{equation}\label{multide}
\Ricc_{\mu}=\left[\begin{smallmatrix}pI&\\&qI\end{smallmatrix}\right]=(2p-q)I+
\left[\begin{smallmatrix}(q-p)I&\\ &2(q-p)I\end{smallmatrix}\right]\in\RR I+\Der(\mu).
\end{equation}
In particular, any bi-invariant complex nilpotent Lie group $(N_{\mu},J)$ (see (\ref{biinv})) admits a
compatible metric which is minimal.

We will now focus on the abelian complex case (see (\ref{abelian})).  It is not hard to see that these
conditions imply that $\Ricc_{\mu}|_{\ngo_1}=\Ric_{\mu}$, and so in this case, to get
$\Ricc_{\mu}|_{\ngo_1}\in\RR I$ is necessary and sufficient that
$$
\la A+F,B\ra=0, \qquad \la A+F,C\ra=0, \qquad ||A||^2=||F||^2.
$$
In order to avoid homothetical changes we will always ask for $||v_1||^2+||v_2||^2=2$, which is equivalent to
$\scalar(\mu)=-1$.

\begin{example} {\rm If we put $A=(s,t)$, $F=(-s,t)$, $B=C=D=E=0$, $s^2+t^2=1$, then the corresponding curve
$\mu_{st}$ of minimal compatible metrics satisfies
$$
\Ric_{\mu_{st}}|_{\ngo_2}=\left[\begin{smallmatrix} s^2&0\\ 0&t^2\end{smallmatrix}\right],
$$
proving that $\{\mu_{st}:s^2+t^2=1,\; 0\leq s\leq \frac{1}{\sqrt{2}}\}$ is a curve of pairwise non-isometric
metrics.  It then follows from Corollary \ref{isoiso} that $(N_{\mu_{st}},J)$ is a curve of pairwise
non-isomorphic abelian complex nilpotent Lie groups.  Recall that $\mu_{st}\simeq\hg_3\oplus\hg_3$ for all $0<s$
and $\mu_{01}\simeq\hg_3\oplus\RR^3$.  }
\end{example}

\begin{example} {\rm For $A=(s,t)$, $F=(-s,t)$, $B=(\unm,0)=C=-D=E$, $s^2+t^2=\unm$, the curve
$\mu_{st}$ of minimal compatible metrics satisfies
$$
\Ric_{\mu_{st}}|_{\ngo_2}=\left[\begin{smallmatrix} s^2+\unm&0\\ 0&t^2\end{smallmatrix}\right],
$$
which implies that the family  $\{\mu_{st}:s^2+t^2=\unm\}$ is pairwise non-isometric. It is easy to see that for
$t\ne 0$, $\mu_{st}$ is isomorphic to the complex Heisenberg Lie algebra, and hence $(N_{\mu_{st}},J)$ defines a
curve of pairwise non-isomorphic abelian complex structures on the Iwasawa manifold. Since
$j_{\mu_{st}}(Z_2)^2\notin\RR I$, we have that the hermitian manifolds $(N_{\mu_{st}},J,\ip)$ are not modified
H-type (see Appendix \ref{app}).  }
\end{example}

\begin{example} {\rm Consider the abelian complex structures defined by $A=-F$, $E=B$ and $D=-C$.  In this case,
the Hermitian manifolds $(N_{\mu},J,\ip)$ are modified H-type and $\mu$ is always isomorphic to the complex
Heisenberg Lie algebra when $v_1,v_2\ne 0$.  In fact, by assuming for simplicity that $\la v_1,v_2\ra=0$, then
$$
j_{\mu}(Z)^2=-\unm(\la Z,Z_1\ra^2 ||v_1||^2+\la Z,Z_2\ra^2 ||v_2||^2)I, \qquad \forall Z\in\ngo_2.
$$
For $A=(s,0)=-F$, $B=(0,t)=E$, $s^2+t^2=1$, $D=C=0$, the corresponding curve $\mu_{st}$ of minimal compatible
metrics satisfies
$$
\Ric_{\mu_{st}}|_{\ngo_2}=\left[\begin{smallmatrix} s^2&0\\ 0&t^2\end{smallmatrix}\right],
$$
and so the family  $\{\mu_{st}:s^2+t^2=1,\; 0\leq s\leq \frac{1}{\sqrt{2}} \}$ is pairwise non-isometric and the
abelian complex structures $(N_{\mu_{st}},J)$ are pairwise non-isomorphic.  Each modified H-type metric is
compatible with two spheres of abelian complex structures of this type which can be described by
$$
\{ \pm v_1\times v_2:v_i\in\RR^3,\; ||v_1||^2=2s^2,\; ||v_2||^2=2t^2,\; \la v_1,v_2\ra=0\},
$$
where $v_1\times v_2$ denotes the vectorial product, but one can see that these structures are all isomorphic to
$(N_{\mu_{st}},J)$ (compare with \cite{KtsSlm}).  We finally recall that $\mu_{01}\simeq\hg_5\oplus\RR$, and so
$\ip$ is a minimal compatible metric for the abelian complex nilpotent Lie group $(N_{\mu_{01}},J)$. }
\end{example}

We finally give a curve of non-abelian complex structures on the Iwasawa manifold.

\begin{example} {\rm Let $\mu_t$ be the curve defined by {\small
$$
\begin{array}{lll}
\mu(X_1,X_3)=-tsZ_2, && \mu(X_2,X_3)=sZ_1, \\
 \mu(X_1,X_4)=sZ_1, && \mu(X_2,X_4)=s(2-t)Z_2,
\end{array} \qquad s=\sqrt{2+t^2+(2-t)^2}, \quad t\in\RR.
$$}
We then have that $A=F=0$, $C=D=(s,0)$, $B=-tJC$, $E=(2-t)JC$, and so $(N_{\mu_t},J)$ is a non-abelian complex
nilpotent Lie group for all $t\in\RR$ (see (\ref{integrable})).  Moreover, $\Ric_{\mu}|_{\ngo_1}$ is diagonal
and hence $\Ricc_{\mu}|_{\ngo_1}$ is a multiple of the identity, which implies that $\ip$ is a minimal
compatible metric for all $(N_{\mu_t},J)$ (see (\ref{multide})).  It follows from
$$
\Ric_{\mu_{t}}|_{\ngo_2}=\unm\left[\begin{smallmatrix} 2s^2&0\\ 0&s^2(t^2+(2-t)^2)\end{smallmatrix}\right],
$$
that the hermitian manifolds $\{(N_{\mu_t},J,\ip):1\leq t<\infty\}$ are pairwise non-isometric since
$$
s^2(t^2+(2-t)^2)-2s^2=(t^2+(2-t)^2)^2-4
$$
is a strictly increasing non-negative function for $1\leq t$, which vanishes if and only if $t=1$.  We therefore
obtain a curve $\{(N_{\mu_t},J):1\leq t<\infty\}$ of pairwise non-isomorphic non-abelian complex nilpotent Lie
groups.  A natural question is which are the nilpotent Lie groups involved.  We have for all $t$ that
$$
j_{\mu_t}(Z_1)=\left[\begin{smallmatrix} &&0&-s\\ &&-s&0\\ 0&s&&\\ s&0&&\end{smallmatrix}\right],\qquad
j_{\mu_t}(Z_2)=\left[\begin{smallmatrix} &&ts&0\\ &&0&-(2-t)s\\ -ts&0&&\\ 0&(2-t)s&&\end{smallmatrix}\right],
$$
and hence $j_{\mu_t}(Z)$ is non-singular if and only if
$$
-t(2-t)\la Z,Z_1\ra^2-\la Z,Z_2\ra^2\ne 0.
$$
This implies that $\mu_t$ is isomorphic to the complex Heisenberg Lie algebra (i.e. when $j_{\mu_t}(Z)$ is
non-singular for any non-zero $Z\in\ngo_2$) if and only if $1\leq t<2$, providing a curve on the Iwasawa
manifold.  Furthermore, $(N_{\mu_1},J)$ is the bi-invariant complex structure and it can be showed by computing
$j_{\mu_t}(Z)^2$ that $(N_{\mu_t},\ip)$ is not modified H-type for any $1<t$.  We finally note that $\mu_2$ is
isomorphic to the group denoted by $(0,0,0,0,12,14+23)$ in \cite{Slm}, and one can easily see by discarding any
other possibility that actually $\mu_t\simeq\mu_2$ for all $2\leq t<\infty$, which gives rise a curve of
pairwise non-isomorphic structures on such a group.  }
\end{example}

Although it has not been mentioned, most of the curves given in this section have been obtained via the
variational method provided by Theorem \ref{equiv1g}, by using an approach very similar to that in Example
\ref{m26}.

\section{Hypercomplex structures}\label{hyper}

Let $N$ be a real $4n$-dimensional nilpotent Lie group with Lie algebra $\ngo$, whose Lie bracket is denoted by
$\mu :\ngo\times\ngo\mapsto\ngo$.  An invariant {\it hypercomplex} structure on $N$ is defined by a triple $\{
J_1,J_2,J_3\}$ of complex structures on $\ngo$ (see Section \ref{complex}) satisfying the quaternion identities
\begin{equation}\label{quat}
J_i^2=-I,\quad i=1,2,3, \qquad J_1J_2=J_3=-J_2J_1.
\end{equation}
An inner product $\ip$ on $\ngo$ is said to be {\it compatible} with $\{ J_1,J_2,J_3\}$, also called an {\it
hyper-hermitian metric}, if
\begin{equation}\label{ortconhyp}
\la J_iX,J_iY\ra=\la X,Y\ra \qquad \forall X,Y\in\ngo, \; i=1,2,3.
\end{equation}
Two hypercomplex nilpotent Lie groups $(N,\{ J_1,J_2,J_3\})$ and $(N',\{ J_1',J_2',J_3'\})$ are said to be {\it
isomorphic} if there exists an automorphism $\vp:\ngo'\mapsto\ngo$ such that
$$
\vp J_i'\vp^{-1}=J_i, \quad i=1,2,3.
$$
For $\gamma=\{ J_1,J_2,J_3\}$ we therefore have that
$$
G_{\gamma}=\Gl(n,\HH)=\{ \vp\in\Gl(4n):\vp J_i=J_i\vp,\; i=1,2,3\}, \qquad K_{\gamma}=\Spe(n),
$$
and the Cartan decomposition of
$$
\ggo_{\gamma}=\glg(n,\HH)=\{ A\in\glg(4n):AJ_i=J_iA,\; i=1,2,3\}
$$
is given by
$$
\glg(n,\HH)=\spg(n)\oplus\pg_{\gamma}, \qquad \pg_{\gamma}=\{ A\in\pg:AJ_i=J_iA,\; i=1,2,3\}.
$$
The invariant Ricci operator for a compatible metric $\ip\in\cca$ is then given by
$$
\Ricg_{\ip}=\unc(\Ric_{\ip}-J_1\Ric_{\ip}J_1-J_2\Ric_{\ip}J_2-J_3\Ric_{\ip}J_3),
$$
and hence condition $\Ricg_{\ip}=0$ can never holds since $\tr{\Ricg_{\ip}}=\scalar(\ip)<0$ for a non-abelian
nilpotent Lie group.

We now fix $\{ J_1,J_2,J_3,\ip\}$ satisfying (\ref{quat}) and (\ref{ortconhyp}) and consider the algebraic
subvariety $\nca_h:=\nca_{\gamma}\subset\nca$ given by
$$
\nca_h=\{\mu\in\nca: J_i \;\mbox{is integrable for}\; i=1,2,3\},
$$
that is, those nilpotent Lie brackets for which $\{ J_1,J_2,J_3\}$ is an hypercomplex structure on the
corresponding nilpotent Lie group $N_{\mu}$.  Thus (\ref{ideg}) identifies each $\mu\in\nca_h$ with the
hyper-hermitian manifold $(N_{\mu},\{ J_1,J_2,J_3\},\ip)$.

\subsection{Hypercomplex $8$-dimensional nilpotent Lie groups}\label{hyper8}  There are no non-abelian nilpotent
Lie groups of dimension $4$ admitting an hypercomplex structure.  In dimension $8$, hypercomplex nilpotent Lie
groups have been determined by I. Dotti and A. Fino in \cite{DttFin0} and \cite{DttFin2}.  They proved the
following strong restrictions on an $8$-dimensional nilpotent Lie algebra $\ngo$ which admits an hypercomplex
structure: $\ngo$ has to be $2$-step nilpotent, $\dim{\mu(\ngo,\ngo)}\leq 4$, there exists a decomposition
$\ngo=\ngo_1\oplus\ngo_2$ such that $\dim{\ngo_i}=4$, $\ngo_i$ is $\{ J_1,J_2,J_3\}$-invariant and
$\mu(\ngo,\ngo)\subset\ngo_2\subset\zg$, where $\zg$ is the center of $\ngo$.  Thus the Lie bracket of $\ngo$ is
just given by a skew-symmetric bilinear form $\mu:\ngo_1\times\ngo_1\mapsto\ngo_2$, and those for which a fixed
$\{ J_1,J_2,J_3\}$ is integrable are also completely described in \cite{DttFin2} as a $16$-dimensional subspace.
Such a description has been successfully used in \cite{DttFin1} to prove that the associated Obata connections
are always flat.

What shall be studied here are the isomorphism classes of such structures and the existence of minimal
compatible metrics.  As we have seen in the previous sections, these two problems are intimately related, and
the hypercomplex case will not be an exception.

For $\ngo_1=\RR^4$ and $\ngo_2=\RR^4$ consider the vector space $W=\Lambda^2\ngo_1^*\otimes\ngo_2$ of all
skew-symmetric bilinear maps $\mu:\ngo_1\times\ngo_1\mapsto\ngo_2$.  Any $8$-dimensional $2$-step nilpotent Lie
algebra with $\dim{\mu(\ngo,\ngo)}\leq 4$ can be modelled in this way.  Fix basis $\{ X_1,X_2,X_3,X_4\}$ and $\{
Z_1,Z_2,Z_3,Z_4\}$ of $\ngo_1$ and $\ngo_2$, respectively.  Each element in $W$ will be denoted as
$\mu=\mu(a_1,...,a_4,...,f_1,...,f_4)$, where {\small
$$
\begin{array}{ll}
\mu(X_1,X_2)=a_1Z_1+a_2Z_2+a_3Z_3+a_4Z_4, & \mu(X_2,X_3)=d_1Z_1+d_2Z_2+d_3Z_3+d_4Z_4,\\
\mu(X_1,X_3)=b_1Z_1+b_2Z_2+b_3Z_3+b_4Z_4, & \mu(X_2,X_4)=e_1Z_1+e_2Z_2+e_3Z_3+e_4Z_4, \\
\mu(X_1,X_4)=c_1Z_1+c_2Z_2+c_3Z_3+c_4Z_4, & \mu(X_3,X_4)=f_1Z_1+f_2Z_2+f_3Z_3+f_4Z_4.
\end{array}
$$ }
The compatible metric will be $\la X_i,X_j\ra=\la Z_i,Z_j\ra=\delta_{ij}$ and the hypercomplex structure will
always act on $\ngo_i$ by
$$
J_1=\left[\begin{smallmatrix} 0&-1&&\\ 1&0&&\\ &&0&-1\\ &&1&0\end{smallmatrix}\right], \quad
J_2=\left[\begin{smallmatrix} &&-1&0\\ &&0&1\\ 1&0&&\\ 0&-1&&\end{smallmatrix}\right],\quad
J_3=\left[\begin{smallmatrix} &&&-1\\ &&-1&\\ &1&&\\ 1&&&\end{smallmatrix}\right].
$$
If $A=(a_1,...,a_4)$, ..., $F=(f_1,...,f_4)$, then it is easy to prove that $J_i$ is integrable for all
$i=1,2,3$ on $N_{\mu}$ (or $(N_{\mu},\{ J_1,J_2,J_3\})$ is a hypercomplex nilpotent Lie group) if and only if
$$
E=B+J_1D+J_1C, \qquad D=-C-J_2A-J_2F, \qquad F=-A-J_3B+J_3E.
$$
If we define $T:=D+C$, then the above conditions are equivalent to
\begin{equation}\label{integrable3}
D=-C+T, \qquad E=B+J_1T, \qquad F=-A+J_2T.
\end{equation}
In order to use a notation as similar as possible to \cite{DttFin1, DttFin2}, we should put
$T=(t_3,t_2,-t_1,t_4)$.  It is easy to check that $(N_{\mu},\{ J_1,J_2,J_3\})$ is {\it abelian} (i.e.
$\mu(J_i\cdot,J_i\cdot)=\mu$, $i=1,2,3$) if and only if $T=0$.  We note that $\dim{W}=24$, and so condition
(\ref{integrable3}) determine a $\Gl(1,\HH)\times\Gl(1,\HH)$-invariant linear subspace $W_h$ of $W$ of dimension
$16$, and a $12$-dimensional subspace $W_{ah}$ if we ask in addition abelian.

For each $\mu\in W$, the Ricci operator of $(N_{\mu},\{ J_1,J_2,J_3\},\ip)$ (see identification (\ref{ideg}))
restricted to $\ngo_1$, $\Ric_{\mu}|_{\ngo_1}$, is given by formula (\ref{ricciw}), and
$$
\Ric_{\mu}|_{\ngo_2}=\unm[\la v_i,v_j\ra], \quad 1\leq i,j\leq 4, \quad v_i:=(a_i,b_i,c_i,d_i,e_i,f_i).
$$
Since the only symmetric transformations of $\ngo_i=\RR^4$ commuting with all the $J_i's$ are the multiplies of
the identity, we obtain that the invariant Ricci operator satisfies $\Ricg_{\mu}|_{\ngo_i}\in\RR I$ for any
$\mu\in W$.  By arguing as in (\ref{multide}), one obtains that any $\mu\in W$ is minimal.  We summarize in the
following proposition the consequences of this fact.

\begin{proposition}
{\rm (i)} Every hypercomplex $8$-dimensional nilpotent Lie group admits a minimal compatible metric, which is
actually its unique compatible left invariant metric up to isometry and scaling.

\no {\rm (ii)} Two hypercomplex $8$-dimensional nilpotent Lie groups $(N_{\mu},\{ J_1,J_2,J_3\})$ and
$(N_{\lambda},\{ J_1,J_2,J_3\})$ are isomorphic if and only if $\mu$ and $c\lambda$ lie in the same
$\Spe(1)\times\Spe(1)$-orbit for some non-zero $c\in\RR$.

\no {\rm (iii)} The moduli space of all $8$-dimensional hypercomplex nilpotent Lie groups up to isomorphism is
parameterized by
$$
\PP W_{h}/\Spe(1)\times\Spe(1).
$$
The representation $W_h$ is equivalent to $(\sug(2)\otimes\RR^4)\oplus\RR^4$, where $\sug(2)$ is the adjoint
representation and $\RR^4$ is the standard representation of $\SU(2)=\Spe(1)$ viewed as real.  Since the
isotropy of an element in general position is finite, the dimension of this quotient is $15-6=9$.

\no {\rm (iv)} The moduli space of all $8$-dimensional abelian hypercomplex nilpotent Lie groups up to
isomorphism is parameterized by
$$
\PP W_{ah}/\Spe(1)\times\Spe(1).
$$
The representation $W_{ah}$ is equivalent to $\sug(2)\otimes\RR^4$, and since the isotropy of an element in
general position is again finite, the dimension of this quotient is $11-6=5$ (see {\rm Example \ref{5g3}} for
explicit $5$-dimensional families).
\end{proposition}

The proofs of the above results follow from Corollary \ref{isoiso}.  It is easy to see that $j_{\mu}(Z)^2\in\RR
I$ for all $Z\in\ngo_2$ if and only if $T=0$, and so only abelian hypercomplex structures can admit (and does)
modified H-type compatible metrics.  Recall that we are considering here a weaker modified H-type condition by
allowing $c(Z)=0$ (see Appendix \ref{app}).

We will give now explicit continuous families of hypercomplex structures on some particular nilpotent Lie
groups.  Let $\ggo_1$, $\ggo_2$ and $\ggo_3$ denote the $8$-dimensional Lie algebras obtained as the direct sum
of an abelian factor and the following H-type Lie algebras: the $5$-dimensional Heisenberg Lie algebra, the
$6$-dimensional complex Heisenberg Lie algebra and the $7$-dimensional quaternionic Heisenberg Lie algebra.  In
order to avoid homothetical changes we will always ask for $||v_1||^2+...+||v_4||^2=2$, which is equivalent to
$\scalar(\mu)=-1$.

\begin{example} {\rm If we put $T=0$, $A=(0,r,0,0)$, $B=(0,0,s,0)$, $C=(0,0,0,t)$, we have for
each $\mu_{rst}$ that
$$
\Ric_{\mu_{rst}}|_{\ngo_2}=\unm\left[\begin{smallmatrix} 0&&&\\ &r^2&&\\ &&s^2&\\
&&&t^2\end{smallmatrix}\right],
$$
and thus the family
$$
\{ (N_{\mu_{rst}},\{ J_1,J_2,J_3\},\ip):0\leq r\leq s\leq t,\; r^2+s^2+t^2=2\}
$$
of minimal compatible metrics is pairwise non-isometric.  This gives rise then a surface of pairwise
non-isomorphic abelian hypercomplex nilpotent Lie groups (see Corollary \ref{isoiso}).  If $0<r$ then
$\mu_{rst}\simeq\ggo_3$, for $r=0<s$ we get a curve on $\ggo_2$ and for $r=s=0$, $t=1$, a single structure on
$\ggo_1$.   }
\end{example}

\begin{example}\label{5g3} {\rm We now set $T=0$ and choose $A,B,C$ such that {\small
$$
v_1=(\frac{1}{\sqrt{2}},0,0,0,0,-\frac{1}{\sqrt{2}}), \; v_2=(0,\sqrt{\frac{3}{8}},0,0,\sqrt{\frac{3}{8}},0), \;
||v_3||^2+||v_4||^2=\frac{1}{4}, \; ||v_3||^2>||v_4||^2.
$$ }
Assume that two of such elements $\lambda=\mu(v_3,v_4)$ and $\lambda'=\mu(v_3',v_4')$ are in the same
$\Spe(1)\times\Spe(1)$-orbit, say $\lambda'=\vp.\lambda$ for some $\vp=(\vp_1,\vp_2)\in\Spe(1)\times\Spe(1)$.
Recall that
$$
j_{\lambda'}(Z)=\vp_1j_{\lambda}(\vp_2^{-1}Z)\vp_1^{-1}, \qquad \forall Z\in\ngo_2,
$$
(see Appendix \ref{app}) and $\la v_i,v_j\ra=-\unm\tr{j_{\mu}(Z_i)j_{\mu}(Z_j)}$, $1\leq i,j\leq 4$.  It then
follows from
$$
||v_1||>||v_2||>||v_3||>||v_4||, \qquad ||v_1'||>||v_2'||>||v_3'||>||v_4'||
$$
that $j_{\vp.\lambda}(Z_i)=\pm j_{\lambda}(Z_i)$ for all $i=1,...,4$, and hence $v_3'=\pm v_3$ and $v_4'=\pm
v_4$.  Thus we have a family of pairwise non-isomorphic abelian hypercomplex structures on $\ggo_3$ depending on
$5$ parameters (see Corollary \ref{isoiso}).  Analogously, we get a $5$-dimensional family on $\ggo_2$ by
putting $v_1=v_2=0$. }
\end{example}

We are now concerned with explicit continuous families of hypercomplex structures which are non-abelian.

\begin{example} {\rm Consider the family defined by {\small
$$
\begin{array}{lll}
\mu_{rst}(X_1,X_2)=rZ_2, && \mu_{rst}(X_2,X_3)=(1-t)Z_4, \\
\mu_{rst}(X_1,X_3)=sZ_3, && \mu_{rst}(X_2,X_4)=-(1-s)Z_3, \\
 \mu_{rst}(X_1,X_4)=tZ_4, && \mu_{rst}(X_3,X_4)=(1-r)Z_2,
\end{array}
$$ }
which is easily seen to satisfy (\ref{integrable3}) but it is not abelian since $t_1=t_2=t_3=0$ but $t_4=1$.
The Ricci operator on the center is given by
$$
\Ric_{\mu_{rst}}|_{\ngo_2}=\unm\left[\begin{smallmatrix} 0&&&\\ &r^2+(1-r)^2&&\\ &&s^2+(1-s)^2&\\
&&&t^2+(1-t)^2\end{smallmatrix}\right],
$$
and hence the family
$$
\left\{ (N_{\mu_{rst}},\{ J_1,J_2,J_3\},\ip):\unm\leq r\leq s\leq t,\quad r^2+s^2+t^2-r-s-t=-\unm\right\}
$$
is pairwise non-isometric.  This gives rise a surface of pairwise non-isomorphic non-abelian hypercomplex
structures on $\ggo_3$ (see Corollary \ref{isoiso}), since $j_{\mu_{rst}}(Z)$ is invertible for any non-zero
$Z\in\ngo_2$. }
\end{example}

\begin{example}
{\rm Let $\mu_t$ be the curve defined for $0\leq t\leq\frac{1}{\sqrt{3}}$ by {\small
$$
\begin{array}{lll}
\mu_{t}(X_1,X_2)=\sqrt{1-3t^3}Z_1+tZ_2, && \mu_{t}(X_2,X_3)=tZ_4, \\
\mu_{t}(X_1,X_3)=tZ_3, && \mu_{t}(X_2,X_4)=-tZ_3, \\
\mu_{t}(X_1,X_4)=tZ_4, && \mu_{t}(X_3,X_4)=-\sqrt{1-3t^3}Z_1+tZ_2.
\end{array}
$$ }
It is easy to check that $(N_{\mu_t},\{ J_1,J_2,J_3\})$ is always non-abelian hypercomplex (recall that
$t_1=t_2=t_3=0$, $t_4=2t$) and the curve is pairwise non-isomorphic since it follows from
$$
\Ric_{\mu_{t}}|_{\ngo_2}=\left[\begin{smallmatrix} 1-3t^2&&&\\ &t^2&&\\ &&t^2&\\
&&&t^2\end{smallmatrix}\right]
$$
that the curve $(N_{\mu_t},\ip)$ is pairwise non-isometric (see Corollary \ref{isoiso}).  The starting and
ending points are $\mu_0\simeq\ggo_1$ and $\mu_{\frac{1}{\sqrt{3}}}\simeq\ggo_3$, respectively, and
$\mu_t\simeq\ug(2)\oplus\CC^2$ for any $0<t<\frac{1}{\sqrt{3}}$ (see \cite{manus} for further information on
these $2$-step nilpotent Lie algebras constructed via representations of compact Lie groups).  }
\end{example}

\section{Einstein solvmanifolds}\label{einstein}

Our goal in this section is to show that the `moment map' approach proposed in this paper can be also applied to
the study of Einstein solvmanifolds.  After a brief overview of such spaces, we will follow the same path used
to study compatible metrics in the previous sections, but by considering the Ricci operator $\Ric_{\ip}$ itself.
In other words, none geometric structure is considered (or $\gamma=0$).

A {\it solvmanifold} is a solvable Lie group $S$ endowed with a left invariant Riemannian metric, and $S$ is
called {\it standard} if $\ag:=\ngo^{\perp}$ is abelian, where $\ngo=[\sg,\sg]$ and $\sg$ is the Lie algebra of
$S$.  Curiously enough, all known examples of non-compact homogeneous Einstein manifolds are isometric to
standard Einstein solvmanifolds.  These spaces have been extensively studied by J. Heber in \cite{Hbr},
obtaining remarkable structure and uniqueness results.

Let $N$ be a nilpotent Lie group with Lie algebra $\ngo$ of dimension $n$, whose Lie bracket is denoted by $\mu
:\ngo\times\ngo\mapsto\ngo$.  We have in this case $\gamma=0$, thus any inner product is `compatible',
$G_{\gamma}=\Gl(n)$, $K_{\gamma}=\Or(n)$, $\pg_{\gamma}=\pg$, $\nca_{\gamma}=\nca$, $\Ricg=\Ric$ and then
condition $\Ricg=0$ is clearly forbidden for non-abelian $N$.  Moreover, the evolution equation is precisely the
Ricci flow and the corresponding invariant Ricci solitons coincide with the metrics studied in \cite{soliton}.

Given a metric nilpotent Lie algebra $(\ngo,\ip)$, a metric solvable Lie algebra $(\sg=\ag\oplus\ngo,\ip')$ is
called a {\it metric solvable extension} of $(\ngo,\ip)$ if the restrictions of the Lie bracket of $\sg$ and the
inner product $\ip'$ to $\ngo$ coincide with the Lie bracket of $\ngo$ and $\ip$, respectively.  It turns out
that for each $(\ngo,\ip)$ there exists a unique rank-one (i.e. $\dim{\ag}=1$) metric solvable extension of
$(\ngo,\ip)$ which stand a chance of being an Einstein space (see \cite{critical}).  This fact turns the study of rank-one
Einstein solvmanifolds into a problem on nilpotent Lie algebras.  More specifically, a nilpotent Lie algebra
$\ngo$ is the nilradical of a standard Einstein solvmanifold if and only if $\ngo$ admits an inner product $\ip$
such that $\Ric_{\ip}=cI+D$ for some $c\in\RR$ and $D\in\Der(\ngo)$, that is, $\ip$ is minimal.

Let $\nca$ be the variety of all nilpotent Lie algebras of dimension $n$.  Fix an inner product $\ip$ on $\ngo$.
Each $\mu\in\nca$ is then identified via (\ref{ideg}) with the Riemannian manifold $(N_{\mu},\ip)$, but we also
have in this case another identification with a solvmanifold: for each $\mu\in\nca$, there exists a unique
rank-one metric solvable extension $S_{\mu}=(S_{\mu},\ip)$ of $(N_{\mu},\ip)$ standing a chance of being Einstein,
and every $(n+1)$-dimensional rank-one Einstein solvmanifold can be modelled as
$S_{\mu}$ for a suitable $\mu\in\nca$.  We recall that the study of standard solvmanifolds reduces to the
rank-one case (see \cite[4.18]{Hbr}).

The functional $F:\PP \nca\mapsto\RR$ given by $F([\mu])=\tr{\Ric_{\mu}^2}/||\mu||^4$ measures how far is the
metric $\mu$ from being Einstein (see \cite{soliton}).

We obtain in this way in Theorem \ref{equiv1g} the uniqueness up to isometry of Einstein metrics on standard
solvable Lie groups proved in \cite{Hbr}, as well as the variational result given in \cite{critical}
characterizing Einstein solvmanifolds as critical points of a natural curvature functional.

Theorem \ref{equiv2g} gives the relationship between Ricci soliton metrics on nilpotent Lie groups and Einstein
solvmanifolds proved in \cite{soliton}.  Part (i) is a new characterization of these privileged metrics,
claiming that they minimize the square norm of  the Ricci tensor along all left invariant metrics on $N$ of a
given scalar curvature.

Finally, we note that the quotient $\nca\kir\Gl(n)$ defined in Section \ref{qgv} parameterizes rank-one Einstein
solvmanifolds of dimension $(n+1)$ up to isometry, as they are precisely the critical points of
$F:||m||^2:\PP\nca\mapsto\RR$ (see Theorem \ref{equiv1g}).  The decomposition of $\nca\kir\Gl(n)$ in categorical
quotients described in (\ref{decomp}) correspond then to the eigenvalue-types introduced in \cite{Hbr} (see
Section \ref{clasif} for further information).

For explicit examples of Einstein solvmanifolds obtained by this method we refer to \cite{Wll}, where it is
proved that any nilpotent Lie group of dimension $6$ admits a minimal metric (see also \cite{finding} for
dimension $\leq 5$ and a $7$-dimensional curve).

\section{On the classification of invariant geometric structures on nilpotent Lie groups}\label{clasif}

There are only two $4$-dimensional symplectic nilpotent Lie groups.  In \cite{KhkGzMdn}, Y. Khakimdjanov, M.
Goze and  A. Medina classified the $6$-dimensional symplectic nilpotent Lie groups up to isomorphism, and have
appeared several continuous families depending on one and two parameters (see also \cite{Mll} for a
classification in the $\NN$-graded filiform case). The complex case is not different, S. Salamon \cite{Slm} determined
which $6$-dimensional nilpotent Lie algebras admit a complex structure, but the moduli space of all complex
structures up to isomorphism on each one of them is still nebulous.  For the Iwasawa manifold for example, such
a moduli space seem to be particularly rich (see Section \ref{complex}).  The set of isomorphism classes of
$2n$-dimensional complex nilpotent Lie groups contains anyway the isomorphism classes of $n$-dimensional
nilpotent Lie algebras over $\CC$, which is far to be achieved for $n\geq 8$.  All the hypercomplex nilpotent
Lie groups have been recently described in dimension $8$ by I. Dotti and A. Fino \cite{DttFin1} as a family
depending on $16$ parameters.  It is proved in Section \ref{hyper8} that we still need $9$ parameters to
describe the moduli space of $8$-dimensional hypercomplex nilpotent Lie groups up to isomorphism, and $5$ just
for the abelian ones.

All this suggests that the moduli space of isomorphism classes of $n$-dimensional class-$\gamma$ nilpotent Lie
groups is probably a very complicated space for most of the classes of geometric structures, even in low
dimensions.  So that, without any hope of an explicit description of such a moduli, what may be said about it?.
Can we show that it is really unmanageable?.  Can we at least find subspaces which are manifolds or algebraic
varieties and obtain lower bounds for its `dimension'?.

This kind of questions belong to invariant theory.  Given a nilpotent Lie group $N$, the set of all
class-$\gamma$ geometric structures on $\ngo$ is parameterized by a relatively nice space: a closed subset of
the symmetric space $\Gl(n)/G_{\gamma}$.  The isomorphism is however determined by the natural action of
$\Aut(\ngo)$, which is a group in general unknown and `very ugly' from an invariant-theoretic point of view
since it is far from being semisimple or reductive.  We then propose to consider the class-$\gamma$
nilpotent Lie groups of a given dimension all together, by using the variety of nilpotent Lie algebras as we did in the
study of compatible metrics in Section \ref{var}.  The advantage of this unified approach is that the group giving
the isomorphism is the reductive Lie group $G_{\gamma}$; the price to pay is that the space
$$
\nca_{\gamma}=\{\mu\in\nca:{\rm IC}(\gamma,\mu)=0\}
$$
where $G_{\gamma}$ is acting on is really wild.  Fortunately, $\nca_{\gamma}$ is at least a real algebraic
variety, and so the classification problem for such structures may be approached by using the tools from
invariant theory given in Section \ref{git}.

\vs

Let $\ngo$, $\nca$, $\gamma$, $\nca_{\gamma}$, $G_{\gamma}$, $K_{\gamma}$ and $N_{\mu}$ be as in Section
\ref{var}.  We may view $\nca_{\gamma}$ as the variety of all class-$\gamma$ $n$-dimensional nilpotent Lie
groups by identifying each element $\mu\in\nca_{\gamma}$ with a class-$\gamma$ nilpotent Lie group,
\begin{equation}\label{ide1}
\mu\longleftrightarrow   (N_{\mu},\gamma),
\end{equation}
where the geometric structure is defined by left invariant translation.  The action of $G_{\gamma}$ on
$\nca_{\gamma}$ has the following interpretation:  each $\vp\in G_{\gamma}$ determines an isomorphism of
class-$\gamma$ geometric structures
$$
(N_{\vp.\mu},\gamma)\mapsto (N_{\mu},\gamma)
$$
by exponentiating the Lie algebra isomorphism $\vp^{-1}:(\ngo,\vp.\mu)\mapsto(\ngo,\mu)$.  In this way, two
class-$\gamma$ structures $\mu,\lambda$ are isomorphic if and only if they live in the same $G_{\gamma}$-orbit,
and hence the quotient
$$
\nca_{\gamma}/G_{\gamma}
$$
parameterizes the moduli space of all $n$-dimensional class-$\gamma$ nilpotent Lie groups up to isomorphism.
Recall that a nilpotent Lie group $N_{\mu}$ admits a class-$\gamma$ structure if and only if the orbit
$\Gl(n).\mu$ meets the variety $\nca_{\gamma}$, and $(\Gl(n).\mu\cap\nca_{\gamma})/G_{\gamma}$ classifies left
invariant class-$\gamma$ geometric structures on $N_{\mu}$ up to isomorphism.

According to the overview given in Section \ref{qgv}, our first natural question would be which are the closed
orbits, or in other words, what kind of class-$\gamma$ nilpotent Lie groups does the categorical quotient
$\nca_{\gamma}\mum G_{\gamma}$ classify?.

Theorem \ref{RS}, (iii) asserts that an orbit $G_{\gamma}.\mu$ is closed if and only if $m_{\gamma}(\lambda)=0$
for some $\lambda\in G_{\gamma}.\mu$, where $m_{\gamma}$ is the moment map for the action of $G_{\gamma}$ on
$\nca_{\gamma}$.  In view of identification (\ref{ideg}) and the formula $m_{\gamma}=8\Ricg$ (see Proposition
\ref{momric}, (ii)), what we are saying is that the categorical quotient $\nca_{\gamma}\mum G_{\gamma}$
parameterizes precisely those class-$\gamma$ nilpotent Lie groups $(N_{\mu},\gamma)$ which admit a compatible
metric $\ip$ statisfying the strong property $\Ricg_{\ip}=0$.  In the case when $\RR I\subset\ggo_{\gamma}$, we
then have that
$$
\nca_{\gamma}\mum G_{\gamma}=\{ (\RR^{n},\gamma)\},
$$
and in the symplectic case, it follows from Proposition \ref{noherm} that for $\mu\in\nca_s$ the orbit
$\Spe(n,\RR).\mu$ can never be closed, unless $\mu=0$, and hence also
$$
\nca_s\mum\Spe(n,\RR)=\{ (\RR^{2n},\omega)\}.
$$
Moreover, $0\in\overline{\Spe(n,\RR).\mu}$ for any $\mu\in\nca_s$ (see Theorem \ref{RS}, (v)), that is, any
$2n$-dimensional symplectic nilpotent Lie group degenerates to $(\RR^{2n},\omega)$.  We also obtain from Theorem
\ref{retraction} the following topological result.

\begin{proposition}
The moduli space $\nca_s/\Spe(n,\RR)$ is contractible.
\end{proposition}

Recall that all these properties are obviously satisfied when $\RR I\subset\ggo_{\gamma}$.

\vs

Due to the absence of closed orbits, it is natural to consider the wider quotient $\nca_{\gamma}\kir G_{\gamma}$
parameterizing orbits which contain a critical point of the functional square norm of the moment map (see
Section \ref{qgv}).  Recall that the moment map $m_{\gamma}:\nca_{\gamma}\mapsto\pg_{\gamma}$ for the action of
$G_{\gamma}$ on $\nca_{\gamma}$ is given by $m_{\gamma}(\mu)=8\Ricg_{\mu}$, and therefore the following nice
`geometric' characterization of this quotient follows from Theorems \ref{equiv1g}, \ref{equiv2g}.

 \begin{proposition}\label{kirmin}
 $\nca_{\gamma}\kir G_{\gamma}$ classifies
precisely those class-$\gamma$ nilpotent Lie groups admitting a minimal compatible metric.
\end{proposition}

Moreover, every class-$\gamma$ nilpotent Lie group degenerates via the negative gradient flow of
$F=||\riccig||^2$ to one of these special structures admitting a minimal metric (see Theorem \ref{ness2}).

We now describe the decomposition of $\nca_{\gamma}\kir G_{\gamma}$ in categorical quotients given in
(\ref{decomp}).  Notation of Section \ref{qgv} will be constantly used from now on.  The proof of the following
rationality result is given in \cite[Section 4]{Nss} for the general case over $\CC$, but it is seen to be valid
over $\RR$ (see also \cite[Theorem 3.5]{strata} for an alternative proof in this particular case which is
clearly valid over $\RR$).

\begin{theorem}\label{rationalization}\cite{Nss}
Let $[\mu]\in \PP \nca_{\gamma}$ be a critical point of $F$, with $\Ricg_{\mu}=c_{\mu}I+D_{\mu}$ for some
$c_{\mu}\in\RR$ and $D_{\mu}\in\Der(\mu)$.  Then there exists $c>0$ such that the eigenvalues of $cD_{\mu}$ are
all integers prime to each other, say $k_1<...<k_r\in\ZZ$ with multiplicities $d_1,...,d_r\in\NN$.
\end{theorem}

\begin{definition}\label{tipo}
{\rm The data set $\eigen$ in the above theorem is called the {\it type} of the critical point $[\mu]$.  }
\end{definition}

It follows from Theorem \ref{rationalization} that the set of types of critical points is in bijection with the
set of strata $S_{A_1},...,S_{A_s}$ given in Theorem \ref{ness2}, and it will be denoted by
$\{\alpha_1,...,\alpha_s\}$.

Fix a type $ \alpha=\alpha_i=\eigen$.  Since $C_{\la A_{\alpha}\ra}\ne\emptyset$, there exists
$\mu\in\nca_{\gamma}$, $||\mu||=1$,  such that $m_{\gamma}([\mu])=8\Ricg_{\mu}=A_{\alpha}$, where
$$
A_{\alpha}=c_{\alpha}I+D_{\alpha}, \quad D_{\alpha}=\left[\begin{smallmatrix}
k_1I_{d_1}&&\\
&\ddots&\\
&&k_rI_{d_r}
\end{smallmatrix}\right], \quad c_{\alpha}=\left\{\begin{smallmatrix}-\frac{\tr{D_{\alpha}}}{n} & {\rm if}\;
\RR I\not\subset\ggo_{\gamma}, \\ \\ -\frac{\tr{D_{\alpha}}+2}{n} & {\rm if}\; \RR I
\subset\ggo_{\gamma},\end{smallmatrix}\right.
$$
and $I_{d_i}$ denotes the $d_i\times d_i$ identity matrix.  The formula for $c_{\alpha}$ follows from taking
trace and using that $\tr{\Ricg_{\mu}}=0$ in the first case and
$\tr{\Ricg_{\mu}}=\tr{\Ric_{\mu}}=-\unc||\mu||^2$ in the second one. Thus the set $C_{\la A_{\alpha}\ra}$ given
in Theorem \ref{ness2} is precisely the set of critical points $[\mu]$ of $F:\PP V\mapsto\RR$ such that
$8\Ricg_{\mu}$ is conjugate to $A_{\alpha}$, and
$$
G_{\alpha}:=G_{A_\alpha}=G_{\gamma}\cap(\Gl(d_1)\times...\times\Gl(d_r)).
$$

\begin{proposition}
For each type $\alpha=\alpha_i=\eigen$ we have that
$$
V_{\alpha}:=V_{A_\alpha}=\{\mu\in V:D_{\alpha}\in\Der(\mu)\}
$$
and
$$
\tilde{G}_{\alpha}:=\tilde{G}_{A_\alpha}=\left\{ g=(g_1,...,g_r)\in
G_{\alpha}:\prod_{i=1}^{r}(\det{g_i})^{k_i}=\det{g}^{-c_{\alpha}}\right\}.
$$
\end{proposition}

\begin{proof} Since $\alpha$ is a type, there exists
$\mu\in\nca_{\gamma}$, $||\mu||=1$,  such that $8\Ricg_{\mu}=A_{\alpha}=c_{\alpha}I+D_{\alpha}$.  It then
follows from (\ref{ricort}) that
$$
||A_{\alpha}||^2=64\tr(\Ricac_{\mu})^2=64\tr(\Ric_{\mu}\Ricg_{\mu})=64(-\frac{1}{32}c_{\alpha})=-2c_{\alpha}.
$$
This implies that $A_{\alpha}.\mu=\frac{||A||^2}{2}\mu$ if and only if
$A_{\alpha}+\frac{||A||^2}{2}I=D_{\alpha}\in\Der(\mu)$, which proves the equality for $V_{\alpha}$.  The formula
for $\tilde{G}_{\alpha}$ follows from the fact that its Lie algebra is given by
$$
\tilde{\ggo}_{\alpha}=\{ A\in\ggo_{\alpha}:\tr{AD_{\alpha}}=-c_{\alpha}\tr{A}\}
$$
(see Section \ref{qgv}).
\end{proof}

We are finally able to give the decomposition of $\nca_{\gamma}\kir G_{\gamma}$ as a disjoint union of
semi-algebraic varieties
$$
\nca_{\gamma}\kir G_{\gamma}=X_1\cup...\cup X_s,
$$
where each $X_i=(V_{\alpha_i}\cap\nca_{\gamma})\mum G_{\alpha_i}$, a categorical quotient.  This allows us to
approach the classification of invariant geometric structures on nilpotent Lie groups by using invariant
theoretic methods.  By considering each $X_i$ separately, we have for instance that geometric structures of type
$\alpha_i$ are precisely the closed $\tilde{G}_{\alpha_i}$-orbits, and two different closed orbits gives rise to
non-isomorphic structures.  We therefore have that two non-isomorphic structures can always been separated by a
$\tilde{G}_{\alpha_i}$-invariant polynomial on $V_{\alpha_i}$.  Moreover, $X_i$ can be described by using a set of
generators and its relations of the ring of such polynomials $\RR[V_{\alpha_i}\cap\nca_{\gamma}]^{\tilde{G}_{\alpha_i}}$
(see \cite{RchSld}).

We do not know how far is $\nca_{\gamma}\kir G_{\gamma}$ from the whole quotient $\nca_{\gamma}/G_{\gamma}$,
that is, from classifying all the geometric structures of class-$\gamma$ of a given dimension.  In view of
Proposition \ref{kirmin}, the crucial question is how strong is for a class-$\gamma$ structure the property of admitting
a minimal compatible metric.  In low
dimensions, the structures seem to be well-disposed to admit a minimal compatible metric, and the only
obstruction we know at this moment is when $\RR I\not\subset\ggo_{\gamma}$ and the nilpotent Lie algebra is characteristically nilpotent, that is,
$\Der(\ngo)$ is nilpotent.

We now give some examples to illustrate the invariant-theoretic approach proposed in this section.  We first
note that all the examples obtained in the complex (see Section \ref{complex}) and hypercomplex (see Section
\ref{hyper}) cases are of type $(1<2;4,2)$ and $(1<2;4,4)$, respectively.   In Example \ref{heisc} appeared almost-complex
structures of types $(2<3<4;4,2n-6,2)$, $(1<2;4,2)$, $(0<1;2,2)$ and $(0;4)$.  For the types we have obtained in the symplectic
case see Table 1.
{\small
\begin{table}\label{types}
$$
\begin{array}{lc}
\hline \\
\mbox{Example} & \mbox{type} \\ \\
\hline\\
\ref{heis}\; (2n\geq 8)& (2<3<4;3,2n-6,3) \\ \\
\ref{heis}\; (2n=6)& (1<2;3,3) \\ \\
\ref{heis}\; (2n=4)& (3<4<6<7;1,1,1,1) \\ \\
\ref{otra4}& (1<2<3<4;1,1,1,1) \\ \\
\ref{abc}& (1<2;3,3) \\ \\
\ref{m26}& (1<2<3<4<5<6;1,1,1,1,1,1) \\
\\ \hline
\end{array}
$$
\caption{Types of symplectic structures}
\end{table} }

At this moment, we are far to be able of describing the types that will appear in each dimension $n$, and even
from get an estimate for the number $s$ of different types.  It would be interesting to do this at least in low
dimensions.  Recall that in this paper a complete description has been achieved only for symplectic and complex
structures of dimension $4$ (but not for the `almost' versions), and in the $8$-dimensional hypercomplex case.  Another natural
question would be if always $k_1\geq 0$.

\begin{example}\label{12nn}
{\rm Consider the type $\alpha=(1<2;n,n)$ in the symplectic case.  If we put $\ngo=\ngo_1\oplus\ngo_2$ with
$\dim{\ngo_i}=n$, then
$$
V_{\alpha}=\Lambda^2\ngo_1^*\otimes\ngo_2=\Lambda(\RR^n)^*\otimes\RR^n, \qquad
D_{\alpha}=\left[\begin{smallmatrix} I_n &\\ & 2I_n\end{smallmatrix}\right], \qquad c_{\alpha}=-\frac{3}{2}.
$$
Let $\{ X_1,...,X_n\}$ and $\{ Z_1,...,Z_n\}$ be basis of $\ngo_1$ and $\ngo_2$ respectively, with dual basis
$\{\alpha_1,...,\alpha_n\}$ and $\{\beta_1,...,\beta_n\}$.  With respect to the inner product $\ip$ which makes
the basis $\{ X_i,Z_i\}$ orthonormal the $2$-form
$$
\omega=\alpha_1\wedge\beta_1+...+\alpha_n\wedge\beta_n
$$
satisfies $\omega=\la\cdot,J\cdot\ra$, where
$$
J=\left[\begin{array}{cc} 0&-I\\ I&0\end{array}\right].
$$
This implies that
$$
G_{\alpha}=\Spe(n,\RR)\cap(\Gl(n)\times\Gl(n))=\{ (g,(g^t)^{-1}):g\in\Gl(n)\}=\Gl(n)
$$
and
$$
\tilde{G}_{\alpha}=\{ (g,(g^t)^{-1}):\det{g}\det{g}^{-2}=1\}=\{ (g,(g^t)^{-1}):g\in\Sl(n)\}=\Sl(n).
$$
A very special feature of this example is that any $\mu\in V_{\alpha}$ satisfies the Jacobi condition, it is
indeed a two-step nilpotent Lie algebra of dimension $2n$ and $\dim{\mu(\ngo,\ngo)}\leq n$. Moreover, any Lie
algebra of this kind can be represented by an element in $V_{\alpha}$.  Thus $\nca_s\cap V_{\alpha}$ (i.e. those
Lie brackets for which $\omega$ is closed) is a $\Gl(n)$-invariant subspace of $V_{\alpha}$, which is easily
seen to be of dimension $\dim{V_{\alpha}}-\binom{n}{3}$.  It follows from the decomposition
$$
V_{\alpha}=\Lambda^2(\RR^n)^*\otimes\RR^n=W_{\epsilon_1+2\epsilon_2}+W_{\epsilon_3}, \qquad W_{\epsilon_3}=
\Lambda^3\RR^n,
$$
in irreducible $\Gl(n)$-modules ($\epsilon_i$: fundamental weights) that $\nca_s\cap
V_{\alpha}=W_{\epsilon_1+2\epsilon_2}$.  We then obtain that the categorical quotient corresponding to the type
$(1<2;n,n)$ parameterizing symplectic nilpotent Lie groups up to isomorphism which admits a minimal compatible
metric is given by
$$
W_{\epsilon_1+2\epsilon_2}\mum\Sl(n).
$$
The description of this quotient is intimately related with the study of the ring
$\RR[W_{\epsilon_1+2\epsilon_2}]^{\Sl(n)}$ of invariant polynomials.  This is a wide open problem for arbitrary
$n$ in invariant theory, even over the complex numbers, and shows that an explicit classification of symplectic
structures on nilpotent Lie groups would not be feasible.  It can be showed however that
$$
\dim{W_{\epsilon_1+2\epsilon_2}\mum\Sl(n)}\geq\frac{1}{3}(n^3-3n^2-n), \qquad \forall\; n\geq 4.
$$
}
\end{example}

\begin{example}\label{12nnc}
{\rm We now consider the type $\alpha=(1<2;2k,2m)$ in the complex case, $n=k+m$.  If we put
$\ngo=\ngo_1\oplus\ngo_2$ with $\dim{\ngo_1}=2k$, $\dim{\ngo_2}=2m$, then
$$
V_{\alpha}=\Lambda^2\ngo_1^*\otimes\ngo_2=\Lambda(\RR^{2k})^*\otimes\RR^{2m}, \qquad
D_{\alpha}=\left[\begin{smallmatrix} I_{2k} &\\ & 2I_{2m}\end{smallmatrix}\right], \quad
c_{\alpha}=-\frac{k+2m+2}{k+m}.
$$
Let $\{ X_1,...,X_{2k}\}$ and $\{ Z_1,...,Z_{2m}\}$ be basis of $\ngo_1$ and $\ngo_2$ respectively, and consider
$J$ acting on each $\ngo_i$ by
$$
\left[\begin{smallmatrix} 0&-1&&&\\ 1&0&&&\\ &&\ddots&&\\ &&&0&-1\\ &&&1&0
\end{smallmatrix}\right]
$$
and the compatible inner product $\ip$ which makes the basis $\{ X_i,Z_i\}$ orthonormal.  This implies that
$$
G_{\alpha}=\Gl(n,\CC)\cap(\Gl(2k)\times\Gl(2m))=\Gl(k,\CC)\times\Gl(m,\CC)
$$
and
$$
\tilde{G}_{\alpha}=\{ (g_1,g_2):\det{g_1}\det{g_2}^{2}=(\det{g_1}\det{g_2})^{-c_{\alpha}}\},
$$
but we can assume that $\tilde{G}_{\alpha}=\Sl(k,\CC)\times\Sl(m,\CC)$ by discarding the elements which act
trivially.  Any $\mu\in V_{\alpha}$ is again a two-step nilpotent Lie algebra of dimension $2n$ and
$\dim{\mu(\ngo,\ngo)}\leq 2m$.  Thus $\nca_c\cap V_{\alpha}$ (i.e. those Lie brackets for which $J$ is
integrable) is a $\Gl(k,\CC)\times\Gl(m,\CC)$-invariant subspace of $V_{\alpha}$, as well as the subspaces
$\nca_{ac}$ and $\nca_{bc}$ of abelian and bi-invariant complex structures, respectively.  One can easily see
that the dimensions of these subspaces are given by
$$
\dim{\nca_c}=k(3k-1)m, \qquad \dim{\nca_{ac}}=2k^2m, \qquad \dim{\nca_{bc}}=k(k-1)m,
$$
and moreover, $\nca_c=\nca_{ac}\oplus\nca_{bc}$.  The decomposition of $V_{\alpha}$ in irreducible
$\Sl(k,\CC)\times\Sl(m,\CC)$-modules is
$$
V_{\alpha}=\Lambda^2(\RR^{2k})^*\otimes\RR^{2m}=(\ug(k)\otimes\RR^{2m})\oplus (W\otimes\RR^{2m})\oplus
(W\otimes\RR^{2m}),
$$
where the action of $\Sl(k,\CC)$ on the skew-hermitian matrices $\ug(k)$ is given by $A\mapsto g^*Ag$,
$\Lambda^2\CC^k=W\oplus W$ as a real representation of $\Sl(k,\CC)$ and $\RR^{2m}$ is the standard
representation of $\Sl(m,\CC)$ on $\CC^m$ viewed as real.  It follows from a simple dimension argument that
$\nca_{ac}=\ug(k)\otimes\RR^{2m}$ and $\nca_{bc}=W\otimes\RR^{2m}$.

We then obtain for the type $(1<2;2k,2m)$, that the categorical quotients parameterizing abelian complex and
bi-invariant complex nilpotent Lie groups up to isomorphism which admits a minimal compatible metric are
respectively given by
$$
\ug(k)\otimes\RR^{2m}\mum\Sl(k,\CC)\times\Sl(m,\CC), \qquad W\otimes\RR^{2m}\mum\Sl(k,\CC)\times\Sl(m,\CC),
$$
and for complex structures
$$
(\ug(k)\otimes\RR^{2m})\oplus (W\otimes\RR^{2m})\mum\Sl(k,\CC)\times\Sl(m,\CC).
$$
The description of these quotients also lead to open problems in invariant theory, except for $m=1$ and maybe for
other small values of $k$ and $m$. }
\end{example}

\begin{example}\label{12hyper}
{\rm
We can argue analogously to the above example for the type given by $(1<2,4k,4m)$ in the abelian hypercomplex case, and
obtain that the corresponding categorical quotient is given by
$$
\spg(k)\otimes\RR^{4m}\mum\Sl(k,\HH)\times\Sl(m,\HH).
$$
Recall that the case $k=m=1$ has been studied in Section \ref{hyper8}, and $\Sl(1,\HH)=\Spe(1)$.  For a type of
the form $\alpha=(k_1<...<k_r;4,...,4)$ in the hypercomplex case, we will always have that $\tilde{G}_{\alpha}$
is compact as it is a closed subgroup of $\Spe(1)\times...\times\Spe(1)$ ($r$ times).  This implies that any
$\mu\in V_{\alpha}\cap\nca_{\gamma}$ admits a minimal compatible metric, which is actually its unique compatible
metric up to isometry and scaling.  Also, the categorical quotient $V_{\alpha}\cap\nca_{\gamma}\mum
\tilde{G}_{\alpha}$ coincides with the whole quotient $V_{\alpha}\cap\nca_{\gamma}/\tilde{G}_{\alpha}$ since any
orbit is closed.  All this makes intriguing enough the study of these particular types.  }
\end{example}

\section{Appendix}\label{app}

We briefly recall in this appendix some features of Riemannian geometry of left invariant metrics on nilpotent
Lie groups.

Consider the vector space $\sym(\ngo)$ of symmetric real valued bilinear forms on $\ngo$, and
$\pca\subset\sym(\ngo)$ the open convex cone of the positive definite ones (inner products), which is naturally
identified with the space of all left invariant Riemannian metrics on $N$.  Every $\ip\in\pca$ induces a natural
inner product $g_{\ip}$ on $\sym(\ngo)$ given by $g_{\ip}(\alpha,\beta)=\tr{A_{\alpha}A_{\beta}}$ for all
$\alpha,\beta\in\sym(\ngo)$, where $\alpha(X,Y)=g(A_{\alpha}X,Y)$.  We endow $\pca$ with the Riemannian metric
$g$ given by $g_{\ip}$ on the tangent space $\tang_{\ip}\pca=\sym(\ngo)$ for any $\ip\in\pca$.  Thus $(\pca,g)$
is isometric to the symmetric space $\Gl(n)/\Or(n)$.  E. Wilson proved that $(N,\ip)$ and $(N,\ip')$ are
isometric if and only if $\ip'=\vp.\ip:=\la\vp^{-1}\cdot,\vp^{-1}\cdot\ra$ for some $\vp\in\Aut(\ngo)$ (see the
proof of \cite[Theorem 3]{Wls}).  Therefore, although the Lie bracket $\mu$ does not play any role in the
definition of a compatible metric, it is crucial in the study of the moduli space of compatible metrics on
$(N,\gamma)$ up to isometry.

The Ricci curvature tensor $\ricci_{\ip}$ and the Ricci operator $\Ric_{\ip}$ of $(N,\ip)$ are given by (see
\cite[7.39]{Bss}),
\begin{equation}\label{ricci}
\begin{array}{rl}
\ricci_{\ip}(X,Y)=\la\Ric_{\ip}X,Y\ra=&-\unm\displaystyle{\sum_{ij}}\la\mu(X,X_i),X_j\ra\la\mu(Y,X_i),X_j\ra \\
&+\unc\displaystyle{\sum_{ij}}\la\mu(X_i,X_j),X\ra\la\mu(X_i,X_j),Y\ra,
\end{array}
\end{equation}
for all $X,Y\in\ngo$, where $\{ X_1,...,X_n\}$ is any orthonormal basis of $(\ngo,\ip)$.  Notice that always
$\scalar(N,\ip)<0$, unless $N$ is abelian.  It is proved in \cite{Jns} that the gradient of the scalar curvature
functional $\scalar:\pca\mapsto\RR$ is given by
\begin{equation}\label{gradsc}
\grad(\scalar)_{\ip}=-\ricci_{\ip},
\end{equation}
and hence it follows from the properties of $\pca$ described above that
\begin{equation}\label{ricort}
\tr{\Ric_{\ip}D}=0, \qquad\forall\; {\rm symmetric}\; D\in\Der(\ngo),
\end{equation}
where $\Der(\ngo)$ is the Lie algebra of derivations of $\ngo$ (see for instance \cite[(2)]{critical} for a
proof of this fact).

Assume now that $\ngo$ is $2$-step nilpotent, and let $\ip$ an inner product on $\ngo$.  Consider the orthogonal
decomposition $\ngo=\ngo_1\oplus\ngo_2$, where $\ngo_2$ is the center of $\ngo$.  Thus the Lie bracket of $\ngo$
can be viewed as a skew-symmetric bilinear map $\mu:\ngo_1\times\ngo_1\mapsto\ngo_2$.  For each $Z\in\ngo_2$ we
define $j_{\mu}(Z):\ngo_1\mapsto\ngo_1$ by
$$
\la j_{\mu}(Z)X,Y\ra=\la\mu(X,Y),Z\ra, \qquad X,Y\in\ngo_1.
$$
$(N,\ip)$ is said to be a {\it modified H-type} Lie group if for any non-zero $Z\in\ngo_2$
$$
j_{\mu}(Z)^2=c(Z)I \qquad \mbox{for some}\; c(Z)<0,
$$
and it is called {\it H-type} when $c(Z)=-\la Z,Z\ra$ for all $Z\in\ngo_2$.  These metrics, introduced by A.
Kaplan, play a remarkable role in the study of Riemannian geometry on nilpotent and solvable Lie groups (see for
instance \cite{BrnTrcVnh} for further information and \cite{modified} for the `modified' case).

If $\mu'=\vp.\mu$ for some $\vp=(\vp_1,\vp_2)\in\Gl(\ngo_1)\times\Gl(\ngo_2)$, then it is easy to see that
$$
j_{\mu'}(Z)=\vp_1j_{\mu}(\vp_2^tZ)\vp_1^t, \qquad \forall Z\in\ngo_2.
$$

\end{document}